\documentclass[12pt,oneside,letterpaper]{amsart}

\usepackage{amsmath,amssymb,amsfonts,amsthm}
\usepackage{mathrsfs}
\usepackage{palatino,mathpazo}
\usepackage[dvipsnames]{xcolor}
\usepackage{cite}
\usepackage[colorlinks,breaklinks,bookmarks=false,linkcolor=NavyBlue,
 urlcolor=ForestGreen,citecolor=ForestGreen,hyperfootnotes=false]{hyperref}
\usepackage[margin=2.0cm,footskip=24pt]{geometry}

\usepackage{enumitem}
\newtheorem{theorem}{Theorem}[section]
\numberwithin{equation}{section}
\newtheorem{introthm}{Theorem}

\newtheorem{proposition}[theorem]{Proposition}
\newtheorem{lemma}[theorem]{Lemma}
\newtheorem{corollary}[theorem]{Corollary}

\theoremstyle{definition}
\newtheorem{definition}[theorem]{Definition}

\newtheorem{remark}[theorem]{Remark}
\newtheorem{example}[theorem]{Example}

\newcommand{\Ric}{\operatorname{Ric}}

\newcommand{\FS}{\mathrm{FS}}
\newcommand{\Vol}{\operatorname{Vol}}
\newcommand{\CP}{\mathbb{CP}}
\newcommand{\BB}{\mathbb{B}}

\title[Toeplitz coefficients and the Lu--Tian correspondence]{Toeplitz coefficients and the Lu--Tian correspondence on Sasakian manifolds}
\author{Matteo Palmieri}
\address{Dipartimento di Matematica e Informatica, Universit\`a degli Studi di Cagliari, Via Ospedale 72,
Cagliari, 09124, Sardinia, Italy}
\email{matteo.palmieri@unica.it}

\author[Wei-Chuan Shen]{Wei-Chuan Shen}
\address{Laboratoire AGM, UMR CNRS 8088, CY Cergy Paris Universit\'e,
2 avenue Adolphe Chauvin, 95302 Cergy-Pontoise Cedex, France}
\curraddr{Department of Mathematics, Syracuse University, 304H Carnegie Hall, Syracuse,
NY 13244-1150, USA}
\email{wshen14@syr.edu}
\thanks{The first author was partially supported by Fondazione di Sardegna through the project ProBiki, and by INdAM--GNSAGA (Gruppo Nazionale per le Strutture Algebriche, Geometriche e loro
Applicazioni). The second author was partially supported by the ANR funded project Quantification des
vari\'et\'es de Caract\`eres comme Mod\`ele pour le chaos quantique (ANR-23-CE40-0021-01).}
\date{}

\theoremstyle{definition}

\begin{document}
\begin{abstract}
Let $X$ be a compact embeddable strictly pseudoconvex CR manifold with Reeb field $T$ and Szeg\H o
projection $\Pi$. Herrmann, Hsiao, Marinescu and the second author showed, for any volume form, that the
semi-classical cut-offs $\chi(k^{-1}\Pi P\Pi)$ of first-order elliptic self-adjoint Toeplitz operators have
full asymptotic expansions. For $P=-\sqrt{-1}T$ and $T$-invariant volume forms, Chang, Herrmann and Hsiao
showed that the diagonal coefficients are $a_{j}\int\chi(t)t^{n-j}dt$ with local invariants $a_{j}$,
computed $a_{0}$ and $a_{1}$, and used normal forms in the Sasakian case. We add one step: on a compact
Sasakian manifold, regular or not, $a_{j}=a_{j}(g^{T})/(2\pi^{n+1})$, where $a_{j}(\cdot)$ is the universal
Bergman kernel polynomial and $g^{T}$ the transverse K\"ahler metric; no line bundle on $X$ is needed.
This gives an exact correspondence: the logarithmic term of the Szeg\H o kernel vanishes to infinite order
over an open set if and only if $a_{j}(g^{T})=0$ there for all $j\ge n+1$, which for circle bundles is a
local form of the theorem of Lu and Tian with its converse. For $(n,0)$-forms on complete Sasakian
manifolds, a cut-off by the self-adjoint generator of the Reeb flow, which avoids the self-adjointness
question for the Toeplitz operator, gives the same results with twisted coefficients.
\end{abstract}
\maketitle

{\small\setcounter{tocdepth}{1}\tableofcontents}

\section{Introduction}\label{sec:intro}

Let $M$ be a compact connected complex manifold of dimension $n$, let $(L,h)$ be a positive
Hermitian line bundle on $M$, and let $X$ be the circle bundle of $(L^{*},h^{*})$, the boundary of the
unit disc bundle $D\subset L^{*}$. The Szeg\H o kernel of $X$, acting on functions, has a boundary
expansion of the form $\Phi\rho^{-n-1}+\Psi\log\rho$, which goes back to Boutet de Monvel and
Sj\"ostrand \cite[Thm.~1.5, Cor.~1.7]{BouSj75}. Under the circle action it decomposes into the Bergman
kernels of the powers $L^{k}$, whose diagonal values have the semi-classical expansion
\cite{Ti90,Zel98,Cat99,Lu00,MM07}
\begin{equation}\label{eq:intro:TYZ}
 B^{L}_{k}\sim\Big(\frac{k}{\pi}\Big)^{n}\sum_{j\ge0}a^{L}_{j}\,k^{-j},\qquad a^{L}_{j}=a_{j}(g),
\end{equation}
where $\omega_{g}=\frac{\sqrt{-1}}{2}\partial\overline\partial\varphi$ is the curvature form of $h$ and the
$a_{j}(\cdot)$ are universal polynomials in the curvature and its covariant derivatives
\cite[Thm.~1.1]{Lu00}. Thus $a_{j}(\cdot)$ is defined for every K\"ahler metric, whether or not it is the
curvature form of a line bundle; below it is evaluated on the transverse K\"ahler metric of a Sasakian
manifold. Lu and Tian proved that the vanishing of the logarithmic term of the
Szeg\H o kernel, which they formulate off the diagonal \cite[(2.4)]{LT}, forces $a^{L}_{j}=0$ for
$j\ge n+1$ \cite[Thm.~2.3]{LT}. The converse is recorded in \cite{LMZ19} as an open conjecture,
attributed to Lu. Loi, Mossa and Zuddas proved it when the expansion \eqref{eq:intro:TYZ} is finite in
the exact sense \cite[Thm.~1.1]{LMZ19}, see also \cite{ALZ13,LUZ16}, and Ebenfelt, Xiao and Xu proved
the corresponding statement for the Bergman kernel of $D$ \cite{EXX}. For compact $M$ the converse also
follows by summation: on the diagonal of $D$ from \eqref{eq:intro:TYZ} with its remainder estimates in
$\mathscr{C}^{m}$ (Theorem~\ref{thm:LuTian}), and on $X\times X$ near the diagonal from the
off-diagonal expansion of the Bergman kernels of $L^{k}$ \cite[Thm.~3.1]{BBS08}, whose coefficients are
almost holomorphic \cite[\S3]{BBS08}.

This summation needs the line bundle. The purpose of this paper is to establish the correspondence
where there is none: on compact Sasakian manifolds, including the irregular ones, whose Reeb flow does
not generate a circle action, and, for $(n,0)$-forms, on Sasakian manifolds which are complete but not
compact.

\subsection{Functions on compact Sasakian manifolds}\label{sec:intro:main}
Let $X$ be a compact strictly pseudoconvex CR manifold of dimension $2n+1$ which is CR embeddable,
equivalently, whose Kohn Laplacian has closed range; let $\eta$ be a contact form with positive Levi form,
$T$ its Reeb field, $dV$ a smooth volume form, and $\Pi$ the orthogonal projection of $L^{2}(X,dV)$ onto
the Hardy space $H^{0}_{b}(X)=\{u\in L^{2}(X,dV):\overline\partial_{b}u=0\}$.

\emph{The general expansion.} Herrmann, Hsiao, Marinescu and the second author proved that for every
formally self-adjoint first-order classical pseudodifferential operator $P$ whose principal symbol is
positive on $\eta$, and every $\chi\in\mathscr{C}^{\infty}_{c}((0,\infty))$, the operator $\chi(k^{-1}\Pi P\Pi)$ is a
semi-classical Fourier integral operator with complex phase, whose kernel has a full asymptotic
expansion as $k\to\infty$; the leading term is $\frac{1}{2\pi^{n+1}}\frac{dV_{\eta}}{dV}\chi(\sigma_{P}(\eta)t)\,t^{n}$,
with $dV_{\eta}=\frac{2^{-n}}{n!}\eta\wedge(d\eta)^{n}$ \cite[Thm.~1.1.1, Cor.~1.1.2]{HHMS23}. No condition on $dV$ is
needed. For $P=-\sqrt{-1}T$, which is formally self-adjoint when $dV$ is invariant under $T$, this is the
semi-classical Toeplitz operator of the Reeb field, $\chi(k^{-1}\Pi(-\sqrt{-1}T)\Pi)$.

\emph{The coefficients.} For a volume form invariant under $T$, Chang, Herrmann and Hsiao showed that the
amplitude depends on $\chi$ only through the derivatives $\chi^{(s)}$, with coefficients which do not depend
on $\chi$ \cite[Thm.~1.3]{CHH25}, and that on the diagonal the expansion reads
$\sum_{j}k^{n+1-j}a_{j}(x)\int\chi(t)t^{n-j}dt$, where the $a_{j}$ are determined by the CR structure, the
contact form and the volume form near $x$ \cite[Thm.~1.5, Rem.~1.6]{CHH25}. They gave the first two
explicitly: $a_{0}=\frac{1}{2\pi^{n+1}}\frac{dV_{\eta}}{dV}$, which also follows from \cite{HHMS23}, and $a_{1}$ in
terms of the Tanaka--Webster scalar curvature $R_{\mathrm{scal}}$ and the density of $dV$; for $dV=dV_{\eta}$,
$a_{1}=R_{\mathrm{scal}}/(4\pi^{n+1})$ \cite[Thms.~1.2, 1.5]{CHH25}. When $T$ preserves the CR structure, that is,
when $(X,\eta,T)$ is Sasakian, they used normal forms: the phase can be taken affine in the Reeb
variables, and the amplitude of the form $\alpha_{j}\chi(t)t^{n-j}$ with $\alpha_{j}$ independent of them
\cite[Thm.~1.3]{CHH25}, determined by its values on the diagonal \cite[(3.1.85), (3.1.86)]{CHH25}. For circle
bundles of positive line bundles, $a_{j}=b_{j}\circ\pi/(2\pi)$ with $b_{j}$ the coefficients of the Bergman
kernel expansion of $L^{k}$ \cite[Thm.~4.2.3]{HHMS23}, \cite[Ex.~1.7]{CHH25}.

\emph{One more step.} For functions on a compact Sasakian manifold, with $dV=dV_{\eta}$, one more step
determines all the $a_{j}$, whether or not the Reeb flow generates a circle action, and gives an exact
correspondence for the Szeg\H o kernel. Here $[T,\Pi]=0$ (Lemma~\ref{lem:sa:commute}), so the operator is
$\chi(-\sqrt{-1}T/k)\Pi$.

\begin{introthm}[Toeplitz coefficients and the Szeg\H o symbol]\label{thm:introD}
Let $(X,\eta,T)$ be a compact Sasakian manifold of dimension $2n+1$, let $\Pi$ be its Szeg\H o projection
on functions, defined with the volume form $dV_{\eta}=\frac{2^{-n}}{n!}\eta\wedge(d\eta)^{n}$, and let
$a_{j}(g^{T})$ be the universal coefficients of \eqref{eq:intro:TYZ} evaluated on the transverse
K\"ahler metric $g^{T}$. Then for every $\chi\in\mathscr{C}^{\infty}_{c}((0,\infty))$, uniformly on $X$,
\begin{equation}\label{eq:intro:sasaki}
 \big(\chi(-\sqrt{-1}T/k)\,\Pi\big)(x,x)\sim\frac{1}{2\pi^{n+1}}\sum_{j\ge0}
 \Big(\int_{0}^{\infty}\chi(t)\,t^{n-j}dt\Big)\,a_{j}(g^{T})(x)\,k^{\,n+1-j}.
\end{equation}
Moreover, in the coordinates of Lemma~\ref{lem:sa:BRT}, modulo a smoothing operator,
\[
 \Pi\equiv\int_{0}^{\infty}e^{\sqrt{-1}t\phi(x,y)}s(z,w,t)\,dt,
\]
where $\phi$ is the phase \eqref{cf:eq:BRTphase} and $s\sim\sum_{j}s_{j}t^{n-j}$ is a classical symbol of order $n$ whose terms do not depend on the Reeb
variables; it is unique modulo flat terms, and
\begin{equation}\label{eq:intro:diag}
 s_{j}(x,x)=\frac{a_{j}(g^{T})(x)}{2\pi^{n+1}}\qquad\text{for every }j\ge0 ,
\end{equation}
so that each $s_{j}$ is, modulo functions flat at $z=w$, the almost analytic extension of
$a_{j}(g^{T})/(2\pi^{n+1})$ from the diagonal.
\end{introthm}
The two parts of Theorem~\ref{thm:introD} are linked by an exact Fourier transform in the Reeb variable
(Proposition~\ref{prop:sa:exact}): in the coordinates of Lemma~\ref{lem:sa:BRT} the flow of $T$ is the
translation $\vartheta\mapsto\vartheta+u$, the phase $\phi$ is linear in the Reeb variables and $s$ does not
depend on them, so $\chi(-\sqrt{-1}T/k)$ multiplies the amplitude by $\chi(t/k)$, and no stationary phase
is needed. For phases linear in the Reeb variables, this structure of the expansion is in
\cite[Thm.~1.3]{CHH25}.
The additional step (Section~\ref{sec:coef}) compares the germ of $X$ at a point with a circle bundle over
$\CP^{n}$ which realizes the same transverse K\"ahler potential: by the locality of the coefficients
\cite[Rem.~1.6]{CHH25} the $a_{j}$ agree at the point, and on the circle bundle they are known
\cite[Thm.~4.2.3]{HHMS23}. The line bundle $\mathcal O(1)\to\CP^{n}$ serves only as a comparison model:
no line bundle on $X$, or on a quotient of $X$, is used, and whether $X$ is regular is a global property
which the germ does not see (Remark~\ref{rem:sa:model}). The logarithmic term of $\Pi$ is defined through the decomposition of Boutet de Monvel
and Sj\"ostrand (Definition~\ref{cf:def:log}).

\begin{introthm}[The logarithmic term]\label{thm:introB}
Let $(X,\eta,T)$, $\Pi$ and $a_{j}(g^{T})$ be as in Theorem~\ref{thm:introD}, and let $W\subset X$ be
open. The logarithmic term of $\Pi$ vanishes to infinite order over $W$ if and only if
$a_{j}(g^{T})=0$ on $W$ for every $j\ge n+1$.
\end{introthm}
This is Theorem~\ref{thm:logS}. For compact $X$, $\int_{X}a_{j}(g^{T})\,dV_{\eta}=0$ for $j\ge n+1$, so the
logarithmic term vanishes if and only if each $a_{j}(g^{T})$, $j\ge n+1$, is constant, if and only if
none of them changes sign (Corollaries~\ref{cor:sa:integrals} and~\ref{cor:sa:compact}). For $n=1$ it
vanishes over a connected open set if and only if $g^{T}$ has constant curvature there
(Corollary~\ref{cor:sa:three}); this local statement uses the second logarithmic coefficient, and for
compact $X$ the consequence that the logarithmic term never vanishes when $X$ is irregular also follows
from Hirachi's formula for the first coefficient \cite{Hir93} (Remark~\ref{rem:sa:compare}). For circle
bundles, $a_{j}(g^{T})=a^{L}_{j}\circ\pi$, and Theorem~\ref{thm:introB} contains the theorem of Lu and
Tian and its converse, in local form; the formulation on the diagonal of $D$ used in
\cite{ALZ13,LUZ16} is equivalent to it (Section~\ref{sec:circle}).

\subsection{$(n,0)$-forms on complete Sasakian manifolds}\label{sec:intro:n0}
Let now $(X,\eta,T)$ be a Sasakian manifold of dimension $2n+1$ whose Sasakian metric is complete; $X$ need
not be compact. Let $H^{n,0}_{b}(X)$ be the Hardy space of $(n,0)$-forms, the $(n,0)$-forms $\alpha$ with
$\overline\partial_{b}\alpha=0$ which are square integrable for the pairing
$\langle\alpha,\beta\rangle=\int_{X}2^{-n}(\sqrt{-1})^{n^{2}}\alpha\wedge\overline\beta\wedge\eta$, and let $S$ be the orthogonal
projection onto it. Hsiao, Marinescu and Wang showed that $S$ is a Fourier integral operator with
complex phase on every relatively compact chart \cite[Thm.~1.3]{HMW22}. On a non-compact manifold the
functional calculus of the Toeplitz operator $S(-\sqrt{-1}\mathcal L_{T})S$ would require its
self-adjointness. We therefore use a slightly different semi-classical cut-off, $\chi(k^{-1}A)S$, where
$A=-\sqrt{-1}\mathcal L_{T}$ is the generator of the unitary group of the Reeb flow on all square
integrable $(n,0)$-forms: $A$ is self-adjoint because the Reeb flow is isometric, and it commutes with $S$
(Lemma~\ref{tf:lem:sa}). On a compact manifold the two cut-offs agree. With this cut-off the remaining
steps go through with the results cited above.

\begin{introthm}[$(n,0)$-forms]\label{thm:introC}
Let $(X,\eta,T)$ be a Sasakian manifold of dimension $2n+1$ whose Sasakian metric is complete, and let
$S$ and $A$ be as above. Let $W\Subset X$ be open, and on a BRT chart in $W$ let $K$ be the kernel of $S$
in the frame $e=dz_{1}\wedge\dots\wedge dz_{n}$, with respect to $dV_{e}=|e|^{2}dV_{\eta}$. Then
$K\equiv\int_{0}^{\infty}e^{\sqrt{-1}t\phi}\sigma\,dt$ with $\sigma\sim\sum_{j}\sigma_{j}t^{n-j}$ independent of
the Reeb variables, and for every $\chi\in\mathscr{C}^{\infty}_{c}((0,\infty))$ the kernel of
$\chi(k^{-1}A)S$ equals $\int_{0}^{\infty}e^{\sqrt{-1}t\phi}\chi(t/k)\sigma\,dt+O(k^{-\infty})$ there and is
$O(k^{-\infty})$ off the diagonal of $W\times W$. Moreover
\[
 |e(x)|^{2}\sigma_{j}(x,x)=\frac{a^{E}_{j}(g^{T})(x)}{2\pi^{n+1}},
\]
where $a^{E}_{j}(\cdot)$ is the universal curvature polynomial giving the $j$-th coefficient of the Bergman
kernel expansion of $L^{k}\otimes K_{M}$ for a positive line bundle $L$ over a K\"ahler manifold $M$, here
evaluated on $g^{T}$ (no line bundle on $X$ is involved); and the logarithmic term of $S$ vanishes to
infinite order over $W$ if and only if $a^{E}_{j}(g^{T})=0$ on $W$ for every $j\ge n+1$.
\end{introthm}
This is Proposition~\ref{tf:prop:H}, Theorem~\ref{tf:thm:main}, Theorem~\ref{tf:thm:id} and
Corollary~\ref{tf:cor:LT}. In the regular case the coefficients $a^{E}_{j}$ also govern the Bergman kernel
of the disc bundle \cite[Thm.~1.9]{EXX} (Proposition~\ref{tf:prop:BS}), and for circle bundles over ball
quotients both logarithmic terms vanish (Corollary~\ref{tf:cor:ball}).

\subsection{Known results used, and what is new}\label{sec:intro:used}
The expansion of the semi-classical Toeplitz operator is \cite{HHMS23}; its form on the diagonal, the
locality and the first two coefficients are \cite{CHH25}; the phase and symbol normal forms go back to
Boutet de Monvel and Sj\"ostrand \cite{BouSj75} and Hsiao and Marinescu \cite{HM14,HM17JDG}, and, in the
Sasakian case, to Chang, Herrmann and Hsiao \cite{CHH25}; and the structure of the Szeg\H o projection on
$(n,0)$-forms of a complete manifold is \cite{HMW22}. Since this paper is addressed also to readers from
K\"ahler geometry, Section~\ref{sec:nf} explains the normal forms, and Appendix~\ref{app:nf} gives short local proofs; being local,
they also apply to the complete manifolds of Theorem~\ref{thm:introC}. Our contribution is the additional
step from these results to the identification of all the coefficients (Theorems~\ref{thm:introD}
and~\ref{thm:introC}) and to the correspondence for the logarithmic term on Sasakian manifolds
(Theorems~\ref{thm:introB} and~\ref{thm:introC}), together with the cut-off by $A$ for $(n,0)$-forms.

\subsection{Organization}
The paper has two parts. The first, Sections~\ref{sec:sasaki}--\ref{sec:compact}, concerns functions on
compact Sasakian manifolds. Section~\ref{sec:sasaki} recalls Sasakian manifolds, BRT coordinates and the
known expansions of the Toeplitz operator of the Reeb field. Section~\ref{sec:nf} explains the normal
forms of the Szeg\H o kernel, for a kernel satisfying local hypotheses which hold in both parts; their
proofs are in Appendix~\ref{app:nf}.
Sections~\ref{sec:coef} and~\ref{sec:log} prove Theorems~\ref{thm:introD} and~\ref{thm:introB}.
Section~\ref{sec:circle} treats the regular case, that is, circle bundles and the disc bundle, by
summation over the Fourier components, and relates Theorem~\ref{thm:introB} to the theorem of Lu and
Tian; apart from Corollary~\ref{cf:cor:rigid} it does not depend on
Sections~\ref{sec:nf}--\ref{sec:log}. Section~\ref{sec:compact} gives the consequences for compact
Sasakian manifolds and for dimension three. The second part concerns $(n,0)$-forms:
Section~\ref{tf:sec:setting} introduces the cut-off by $A$, describes the Szeg\H o kernel in a frame,
localizes the cut-off and proves Theorem~\ref{thm:introC}, and Section~\ref{tf:sec:regular} treats the
regular case, with the Bergman kernels of $L^{k}\otimes K_{M}$ and of the disc bundle, and ball
quotients.

\subsection{Notation}\label{sec:intro:notation}
\begin{itemize}[leftmargin=*]
\item On a Sasakian manifold $dV_{\eta}=\frac{2^{-n}}{n!}\eta\wedge(d\eta)^{n}$. For the circle bundle of
      $(L^{*},h^{*})$, $\eta$ is the connection form, $T$ generates the circle action and
      $dV_{\eta}=d\theta\wedge\pi^{*}dV_{g}$ with $dV_{g}=\omega_{g}^{n}/n!$.
\item $\Pi$ is the Szeg\H o projection onto the Hardy space $H^{0}_{b}(X)$ of $L^{2}$ CR functions, for
      $dV_{\eta}$. In Section~\ref{sec:circle}, on circle bundles, it is written $S$, with Fourier
      components $S_{k}$ onto $H^{0}_{b,k}(X)$, as in \cite{LT}. From Section~\ref{tf:sec:setting} on, $S$ is
      the Szeg\H o projection onto the Hardy space $H^{n,0}_{b}(X)$ of $L^{2}$ CR $(n,0)$-forms.
\item A symbol is \emph{Reeb-independent} if, in BRT coordinates, its terms do not depend on the Reeb
      variables $\vartheta_{x}$ and $\vartheta_{y}$, that is, depend on $(z,w)$ only.
\item $B^{L}_{k}$ and $B_{k}$ are the Bergman kernel functions of $L^{k}$ and of $E_{k}=L^{k}\otimes K_{M}$,
      with coefficients $a^{L}_{j}$ and $a^{E}_{j}$ as in \eqref{eq:intro:TYZ}; on circle bundles $S_{k}(x,x)=B^{L}_{k}/(2\pi)$ \eqref{eq:SkBk}.
\item A function near the diagonal is \emph{flat} if it vanishes to infinite order there;
      $A\equiv B$ means that $A-B$ is smooth.
\end{itemize}

\section{Sasakian manifolds and BRT coordinates}\label{sec:sasaki}
\subsection{Sasakian structures and the Szeg\H o projection}\label{ssec:sa:setup}
Let $(X,T^{1,0}X)$ be compact, connected and
strongly pseudoconvex of dimension $2n+1$, and let $T$ be a transversal CR
vector field: $[T,\mathscr{C}^{\infty}(X,T^{1,0}X)]\subset\mathscr{C}^{\infty}(X,T^{1,0}X)$
and $\mathbb CTX=T^{1,0}X\oplus T^{0,1}X\oplus\mathbb CT$. Let $\eta$ be the
contact form with $\eta(T)=1$ and $\eta(T^{1,0}X)=0$; replacing $T$ by $-T$ if
necessary, we assume that the Levi form of $\eta$ is positive definite. Then the flow of $T$
preserves $\eta$, and $(X,\eta,T)$ is a Sasakian manifold with transverse
K\"ahler form $\omega^{T}:=\frac12d\eta$ and transverse metric $g^{T}$. Let
$H^{0}_{b}(X)=\{u\in L^{2}(X,dV_{\eta}):\overline\partial_{b}u=0\}$ be the Hardy space, with $\overline\partial_{b}u$
taken in the sense of distributions, and let $\Pi:L^{2}(X,dV_{\eta})\to H^{0}_{b}(X)$ be the orthogonal
projection, the Szeg\H o projection \cite[(1.1.4)]{HHMS23}. By Lemma~\ref{lem:sa:commute},
\begin{equation}\label{eq:sa:commute}
 [T,\Pi]=0 .
\end{equation}
The circle bundle of a positive line bundle is the case in which $T$ generates a free circle action (Section~\ref{sec:circle}). If $T$ generates a locally free circle action,
$X$ is \emph{quasi-regular}, and otherwise \emph{irregular}. Fix
$\chi\in \mathscr{C}^{\infty}_{c}((0,\infty))$ with $\chi\ge0$, $\chi\not\equiv0$, and put
$m_{\ell}(\chi)=\int_{0}^{\infty}\chi(t)\,t^{\ell}\,dt$, so that $m_{\ell}(\chi)>0$ for every
$\ell\in\mathbb Z$. By \eqref{eq:sa:commute}, the semi-classical Toeplitz operator
of $P=-\sqrt{-1}T$ \cite[Thm.~1.1.1]{HHMS23} is
\begin{equation}\label{eq:sa:chiT}
 \chi\big(k^{-1}\Pi(-\sqrt{-1}T)\Pi\big)=\chi(-\sqrt{-1}T/k)\,\Pi
 =\int_{\mathbb R}K_{k}(u)\,e^{uT}\Pi\,du,\qquad K_{k}(u)=\tfrac{k}{2\pi}\widehat{\chi}(ku).
\end{equation}

\begin{lemma}[Invariance of $\Pi$]\label{lem:sa:commute}
Let $(X,\eta,T)$ be a Sasakian manifold, that is, a strictly pseudoconvex CR manifold with a contact form $\eta$ annihilating $T^{1,0}X$ whose Reeb field $T$ preserves $T^{1,0}X$; for instance the circle bundle $X$ with its connection form. Assume that the flow
$\varphi_{s}$ of $T$ is complete, as it is when $X$ is compact. Then $\varphi_{s}$ preserves $T^{1,0}X$ and $dV_{\eta}$, and $U(s)=\varphi_{s}^{*}$
commutes with $\Pi$ for all $s$. In particular $\Pi(\varphi_{s}x,\varphi_{s}y)=\Pi(x,y)$.
\end{lemma}

\begin{proof}
Since $\eta(T)=1$ and $\iota_{T}d\eta=0$, Cartan's formula gives $\mathcal L_{T}\eta=0$,
hence $\mathcal L_{T}dV_{\eta}=0$; and $T$ preserves $T^{1,0}X$, so $\varphi_{s}$ is a CR
diffeomorphism. Thus $U(s)$ is unitary on $L^{2}(X,dV_{\eta})$, and it maps
$\ker\overline{\partial}_{b}$ into itself, because the pull-back of an $L^{2}$
function annihilated by $\overline{\partial}_{b}$ in the sense of distributions is
again such a function. The same holds for $U(-s)=U(s)^{-1}$, so $U(s)$ maps
the closed subspace $\ker\overline{\partial}_{b}$ onto itself; being unitary, it maps
the orthogonal complement onto itself, and therefore commutes with the
orthogonal projection $\Pi$.
\end{proof}

\begin{remark}
The invariance of $dV_{\eta}$ is essential: for a volume form which is not invariant under
the Reeb flow, $U(s)$ still preserves $\ker\overline{\partial}_{b}$, but in general not its orthogonal
complement, and $\Pi$ need not commute with $T$.
\end{remark}

\subsection{The Sasaki cone and closed range}\label{ssec:sa:cone}
We collect the known results used below, with the short proofs needed in
the present setting.

\begin{lemma}[The Sasaki cone]\label{lem:sa:cone}
There are a torus $\mathbb T$ of CR automorphisms of $X$, whose Lie algebra
$\mathfrak t$ contains $T$, and the open convex cone
$\mathfrak t^{+}=\{T'\in\mathfrak t:\eta(T')>0\}$, such that every $T'\in\mathfrak t^{+}$ is the
Reeb field of the Sasakian structure $\eta_{T'}=\eta/\eta(T')$, with the same CR
structure, and such that the $T'\in\mathfrak t^{+}$ generating a circle action are
dense in $\mathfrak t^{+}$.
\end{lemma}
\begin{proof}
Let $\mathbb T$ be the closure of the flow of $T$ in the compact group of
automorphisms of the Sasakian structure; that a transversal CR $\mathbb R$-action
comes from a torus action in this way is \cite[Thm.~3.5, Cor.~3.9]{HeHsLiR},
cf.\ \cite[(2.12)]{HLM23}. For $T'\in\mathfrak t^{+}$ the flow of $T'$
preserves $\eta$, hence $\iota_{T'}d\eta=-d(\eta(T'))$ by Cartan's formula and $T'(\eta(T'))=0$; these
give $\eta_{T'}(T')=1$ and $\iota_{T'}d\eta_{T'}=0$. Since $T'$ is CR and transversal,
it is the Reeb field of a Sasakian structure; this is the Sasaki cone of
\cite{BoGa08}. The elements of $\mathfrak t$ generating circle subgroups are those
rational with respect to the lattice $\ker(\exp:\mathfrak t\to\mathbb T)$, and they
are dense; cf.\ \cite{Ruk95}.
\end{proof}

\begin{lemma}[Closed range]\label{lem:sa:closed}
The operator $\overline{\partial}_{b}$ on $X$ has closed range in $L^{2}$, for every
smooth positive volume form.
\end{lemma}
\begin{proof}
For $n\ge2$ this holds on every compact strongly pseudoconvex manifold
\cite{Koh86}. For $n=1$, by Lemma~\ref{lem:sa:cone} the CR structure of $X$ is
that of a manifold with a transversal CR locally free circle action, which
is CR embeddable \cite{Lem92}, \cite[Thm.~1.13]{HM17JDG}, see also \cite[Thm.~1.2]{HeHsLi17},
so the range is closed \cite{Koh86}; closed range does not depend on the
volume form.
\end{proof}

\subsection{BRT coordinates}\label{ssec:sa:BRT}
The coordinates of Baouendi--Rothschild--Treves \cite{BRT85} are usually
quoted; for transversal CR circle actions they are used in
\cite[Theorem~6.5]{HM17JDG}, and for $\mathbb R$-actions in
\cite[Definition~2.6]{HLM23}. In our situation the proof is short.

\begin{lemma}[BRT coordinates]\label{lem:sa:BRT}
Every point of $X$ has a neighborhood with coordinates
$(z,\vartheta)\in U\times(-\varepsilon,\varepsilon)$, $U\subset\mathbb C^{n}$ open, in which
$T=\partial_{\vartheta}$ and $T^{1,0}X$ is spanned by
$Z_{j}=\partial_{z_{j}}+\frac{\sqrt{-1}}{2}\varphi_{z_{j}}(z)\partial_{\vartheta}$ for a strictly
plurisubharmonic $\varphi$ on $U$. Then
$\eta=d\vartheta+\operatorname{Im}(\partial\varphi)$ and $d\eta=\sqrt{-1}\partial\overline{\partial}\varphi$.
\end{lemma}
\begin{proof}
Since $T\ne0$, there are coordinates $(y,\vartheta)$ with $T=\partial_{\vartheta}$; let
$\pi(y,\vartheta)=y$. Let $H=\operatorname{Re}(T^{1,0}X\oplus T^{0,1}X)$, with the complex
structure $J$. As $TX=H\oplus\mathbb RT$, $d\pi|_{H}$ is an isomorphism, and as the
flow of $T$ preserves $H$ and $J$, the structure $J_{W}=d\pi\,J\,(d\pi|_{H})^{-1}$ on
the local leaf space $W$ is well defined. A $(1,0)$-field $A$ on $W$ has a
unique lift $\widehat{A}$ in $T^{1,0}X$, which is $T$-invariant; $[\widehat{A},\widehat{B}]$ lies in
$T^{1,0}X$ and is $\pi$-related to $[A,B]$, so $J_{W}$ is integrable, and the
Newlander--Nirenberg theorem \cite{NN57} gives holomorphic coordinates $z$ on $W$. The lift
of $\partial_{z_{j}}$ is $Z_{j}=\partial_{z_{j}}+a_{j}\partial_{\vartheta}$. Since $[\partial_{\vartheta},Z_{j}]$ and $[Z_{j},Z_{k}]$
lie in $T^{1,0}X\cap\mathbb CT=0$, $a_{j}=a_{j}(z)$ and $\alpha=\sum_{j}a_{j}dz_{j}$ satisfies
$\partial\alpha=0$; on a polydisc $\alpha=\partial h$ by the Dolbeault--Grothendieck lemma.
With $h=u+\sqrt{-1}v$ and $\vartheta'=\vartheta-u$ one gets $Z_{j}=\partial_{z_{j}}+\sqrt{-1}v_{z_{j}}\partial_{\vartheta'}$,
which is the claim with $\varphi=2v$. The formulas for $\eta$ and $d\eta$ follow,
and $\varphi$ is strictly plurisubharmonic by strong pseudoconvexity.
\end{proof}

Two
charts differ by $z'=F(z)$ holomorphic, $\vartheta'=\vartheta-\frac12\operatorname{Im}g(z)$ and
$\varphi'=\varphi+\operatorname{Re}g$ with $g$ holomorphic. Since $\operatorname{Re}g$ is pluriharmonic,
the K\"ahler metric with K\"ahler form $\frac{\sqrt{-1}}{2}\partial\overline\partial\varphi$ is the same in both
charts: it is the transverse K\"ahler metric $g^{T}$, and every local curvature invariant of it, such as
$a_{j}(g^{T})$, is a smooth function on $X$, whether or not $X$ is regular. The BRT phase \eqref{cf:eq:BRTphase} and the class $\mathcal C$ of Section~\ref{sec:nf} are
defined in such a chart; they do not depend on the chart (Lemma~\ref{lem:sa:chart}).

\begin{lemma}[Change of chart]\label{lem:sa:chart}
Under a change of chart as above, the phases $\phi$ and $\phi'$ of the two
charts differ by a function that is flat at $z=w$, and a symbol belongs to
the class $\mathcal C$ of one chart if and only if it belongs to that of the
other. For every point $x_{0}$ of a chart there is a change of chart after
which $z(x_{0})=0$ and $\varphi(z)=|z|^{2}+O(|z|^{3})$.
\end{lemma}
\begin{proof}
Regard $z$ and $z'=F(z)$ as two coordinates on the same leaf space. The
function $\varphi(z,\overline{w})+\frac12\big(g(z)+\overline{g(w)}\big)$ is an almost analytic
extension of $\varphi'=\varphi+\operatorname{Re}g$, and two almost analytic extensions differ by
a function flat at $z=w$. Hence, modulo such functions,
\begin{align*}
 \mathcal B'(z,w)&=\mathcal B(z,w)-\sqrt{-1}\Big(\tfrac12g(z)+\tfrac12\overline{g(w)}
 -\tfrac12\operatorname{Re}g(z)-\tfrac12\operatorname{Re}g(w)\Big)\\
 &=\mathcal B(z,w)+\tfrac12\operatorname{Im}g(z)-\tfrac12\operatorname{Im}g(w),
\end{align*}
while $\vartheta'_{x}-\vartheta'_{y}=\vartheta_{x}-\vartheta_{y}-\frac12\operatorname{Im}g(z)+\frac12\operatorname{Im}g(w)$; the sum is
$\phi'=\phi$. A function of $(z,w)$ is a function of $(z',w')$, which gives
the statement on $\mathcal C$. For the last assertion, translate $z$ so that
$z(x_{0})=0$ and let
\[
 P(z)=\varphi(0)+2\sum_{j}\varphi_{z_{j}}(0)z_{j}+\sum_{j,k}\varphi_{z_{j}z_{k}}(0)z_{j}z_{k} .
\]
Taylor's formula gives
$\varphi-\operatorname{Re}P=\sum_{j,k}\varphi_{z_{j}\bar z_{k}}(0)z_{j}\overline{z}_{k}+O(|z|^{3})$. Take $g=-P$
and a linear map $F$ that brings the positive definite Hermitian form on
the right to $|z'|^{2}$.
\end{proof}

\subsection{Known expansions}\label{ssec:sa:known}
We recall the results on which Section~\ref{sec:coef} builds, in the order in which they are used.

\emph{The general expansion.} For any compact embeddable strictly pseudoconvex CR manifold, any smooth
volume form and any formally self-adjoint first-order classical pseudodifferential operator $P$ whose
principal symbol is positive on the contact form, $\chi(k^{-1}\Pi P\Pi)$ is a semi-classical Fourier
integral operator with complex phase, with a full asymptotic expansion and the leading term recalled in
Section~\ref{sec:intro:main} \cite[Thm.~1.1.1, Cor.~1.1.2]{HHMS23}. In the regular case all the coefficients
are identified with those of \eqref{eq:intro:TYZ} \cite[Thm.~4.2.3]{HHMS23}, and for quasi-regular $X$ the
weighted Szeg\H o kernel functions have a smooth expansion with leading coefficient
$\frac{1}{2\pi^{n+1}}\frac{dV_{\eta}}{dV}\int\chi(t)t^{n}dt$ \cite[Thm.~4.3.1]{HHMS23}.

\emph{The coefficients.} For the Reeb field of any contact form and any volume form invariant under it,
the expansion on the diagonal has the form \eqref{eq:sa:diag}, with coefficients which do not depend on
$\chi$ and are locally determined, and the first two are computed, in
\cite[Thms.~1.2, 1.3, 1.5, Rem.~1.6]{CHH25}. When the Reeb field preserves the CR structure, the phase can
be taken affine in the Reeb variables and the amplitude independent of them, and the proof shows that
the amplitude is then determined modulo flat terms by its diagonal values
\cite[Thm.~1.3, (3.1.85), (3.1.86)]{CHH25}.

\emph{Related expansions.} In general $\chi(-\sqrt{-1}T/k)\Pi=\sum_{\lambda}\chi(\lambda/k)\Pi_{\lambda}$, the sum over
the eigenvalues $\lambda$ of $-\sqrt{-1}T$ on $H^{0}_{b}(X)$ \cite[Thms.~4.1, 4.5]{HeHsLiR}; for a locally free
circle action the $\lambda$ are integers, and this is the weighted sum $\sum_{m}\chi(m/k)\Pi_{m}$ of the
Fourier components. Expansions of such weighted Szeg\H o kernels were obtained by Hsiao, Li and
Marinescu \cite[Thm.~1.1]{HLM21} for circle actions and by Herrmann, Hsiao and Li \cite[Thm.~1.1]{HeHsLiR}
for $\mathbb R$-actions, see also \cite{HHL20}, for the tensor powers of a rigid positive CR line bundle; see
\cite[Thm.~5.4, Rem.~1.4]{HLM23} for the form in which they apply to the trivial bundle here. First
coefficients for quasi-regular $X$ are computed in \cite{HeHsLi18}. For quasi-regular $X$, the Bergman
kernels $B_{k}$ of the quotient orbifold behave periodically in $k$ at the orbifold points \cite{DLM06}; the
weighted Bergman kernels $\sum_{i}c_{i}B_{k+i}$ remove this periodicity and have a smooth global expansion
with the coefficients of the smooth case \cite[Thm.~1.7]{RT11_1}, which also follows from
\cite[(5.25)]{DLM06} \cite{DLM11}; the first coefficients are written out in \cite{RT11_1,DLM11,HeHsLi18}. See
also the thesis \cite{Cha26}. In the quasi-regular case the weight is essential: a single $\Pi_{m}(x,x)$
has no expansion that is uniform near the points with non-trivial stabilizer \cite[Thm.~1.1]{HeHsLi18}.

Proposition~\ref{prop:sa:exact} gives the expansion in the present setting, with the Reeb-independent
Szeg\H o symbol as amplitude, and so relates the coefficients to that symbol.

\section{Normal forms of the Szeg\H o kernel}\label{sec:nf}
This section explains the normal forms of the Szeg\H o kernel on which the paper rests: the phase can
be taken affine in the Reeb variable and the symbol independent of it; the symbol is then unique modulo
flat terms and almost sesquiholomorphic, so that it is determined by its values on the diagonal. These
facts are known: the phase and the Reeb-independent symbol for circle bundles are in \cite{HM14}, the
uniqueness argument is that of \cite[\S8]{HM17JDG}, the Reeb-independent amplitude of the semi-classical
Toeplitz operator of a compact Sasakian manifold, with its uniqueness and its determination by the
diagonal values, is in \cite{CHH25}, and the polarization is that of \cite{HS20}; precise references are
given in Remark~\ref{nf:rem:sources}. We state them for a kernel $K$ satisfying the hypotheses
\ref{K:rep}--\ref{K:sym} below, which hold for the Szeg\H o projection on functions of a compact
Sasakian manifold and, in a frame, for the Szeg\H o projection on $(n,0)$-forms of a complete one
(Proposition~\ref{tf:prop:H}). All arguments are local; the proofs are collected in Appendix~\ref{app:nf}.

\subsection{The parametrix of Boutet de Monvel and Sj\"ostrand}\label{nf:sec:param}

The following theorem is due to Boutet de Monvel and Sj\"ostrand
\cite[Thm.~1.5]{BouSj75}; we use it in the form \cite[Thm.~1.9 with $q=0$, Thm.~5.1]{HM17JDG},
where the properties ($\alpha$)--($\gamma$) of the phase are
\cite[(3.2)--(3.4), Thm.~3.4]{HM17JDG}.

\begin{theorem}[Boutet de Monvel--Sj\"ostrand]\label{thm:BdMS}
Modulo a smoothing operator,
\begin{equation}
 \Pi(x,y)=\int_{0}^{\infty}e^{\sqrt{-1}\,t\,\phi(x,y)}\,s(x,y,t)\,dt ,
\end{equation}
where the phase $\phi$ has the following properties.
\begin{itemize}[leftmargin=*]
\item[($\alpha$)] $\operatorname{Im}\phi\ge0$; $\phi=0$ only on the diagonal;
      $\phi(x,y)=-\overline{\phi(y,x)}$; and
      $d_{x}\phi|_{\rm diag}=-d_{y}\phi|_{\rm diag}=\eta$.
\item[($\beta$)] The eikonal condition: $p_{0}(x,\partial_{x}\phi)-f\phi$
      vanishes to infinite order on the diagonal for some smooth $f$, where
      $p_{0}$ is the principal symbol of $\Box_{b}$.
\item[($\gamma$)] At each diagonal point, in coordinates in which the Levi
      form is diagonal, $\phi$ has the Taylor expansion to second order of
      \cite[Thm.~3.4]{HM17JDG}; in particular $\operatorname{Im}\phi(x,y)\ge c|z-w|^{2}$,
      where $z$ and $w$ are the horizontal coordinates of $x$ and $y$ in
      Lemma~\ref{lem:sa:BRT}.
\end{itemize}
Further, $s(x,y,t)\sim\sum_{j\ge0}s_{j}(x,y)\,t^{\,n-j}$ is a classical symbol of
order $n$.
\end{theorem}
In the convention of \cite{HM17JDG}, $\omega_{0}=-\eta$ and
$d_{x}\phi|_{\mathrm{diag}}=-\omega_{0}$. The imaginary part of $\phi$ need not
be positive off the diagonal: along a Reeb orbit the BRT phase
\eqref{cf:eq:BRTphase} is real.

By Lemma~\ref{cf:lem:J}, the term $s_{j}t^{n-j}$ contributes a multiple of
$(-\sqrt{-1}\phi)^{-(n+1-j)}$ for $j\le n$, and, for $j=n+1+p$, a multiple of
$\phi^{p}\log(-\sqrt{-1}\phi)$, modulo smooth functions. So the logarithmic term
of $\Pi$ is carried by $s_{n+1},s_{n+2},\dots$.

The leading coefficient of the symbol is computed in \cite[Thm.~4.5]{HM17JDG};
see also \cite[(1.1.2)--(1.1.5)]{CHH25}.

\subsection{The BRT phase}\label{cf:sec:setting}
Let $D$ be a chart of Lemma~\ref{lem:sa:BRT}, with BRT coordinates $x=(z,\vartheta)$, in which
$T=\partial_{\vartheta}$ and
$\overline{Z}_{l}=\partial_{\bar z_{l}}-\frac{\sqrt{-1}}{2}\varphi_{\bar z_{l}}\partial_{\vartheta}$.
Let $\varphi(z,\overline{w})$ be an almost analytic extension of the potential with
$\overline{\varphi(z,\overline{w})}=\varphi(w,\overline{z})$, and let
\begin{equation}\label{cf:eq:BRTphase}
 \phi(x,y)=\vartheta_{x}-\vartheta_{y}+\mathcal B(z,w),\qquad
 \mathcal B(z,w)=-\sqrt{-1}\big(\varphi(z,\overline{w})-\tfrac12\varphi(z)-\tfrac12\varphi(w)\big)
\end{equation}
be the \emph{BRT phase}. It has the following properties.
\begin{enumerate}[label=\rm(P\arabic*)]
\item $\phi$ is affine in $\vartheta_{x}$ and in $\vartheta_{y}$, and is
      invariant under $(x,y)\mapsto(e^{uT}x,e^{uT}y)$ for small $u$.
\item $\phi(x,y)=-\overline{\phi(y,x)}$, and $\phi(e^{\theta T}x,x)=\theta$ for small $\theta$.
\item $c_{0}|z-w|^{2}\le\delta(z,w):=\operatorname{Im}\mathcal B(z,w)\le C_{0}|z-w|^{2}$
      and $|\mathcal B(z,w)|\le C_{0}|z-w|$.
\item $\overline{Z}_{l,x}\phi=-\sqrt{-1}\,\varepsilon_{l}$, where
      $\varepsilon_{l}(z,w):=\partial_{\bar z_{l}}\varphi(z,\overline{w})$ is flat.
\end{enumerate}
We put $c:=\vartheta_{y}-\mathcal B$, so that $\phi=\vartheta_{x}-c$ and $\operatorname{Im}c=-\delta$.
The upper bound in (P3) is explained in Remark~\ref{cf:rem:upper}.

\subsection{Conventions, inputs and hypotheses}\label{nf:sec:hyp}
\paragraph{Conventions.} A function near the diagonal is \emph{flat} if it is
$O(|x-y|^{N})$ for every $N$, locally uniformly; for functions of $(z,w)$
this means $O(|z-w|^{N})$. A classical symbol $a\sim\sum_{j}a_{j}t^{m-j}$ is
flat if all its terms are. Let $\kappa\in \mathscr{C}^{\infty}(\mathbb R)$ vanish on
$(-\infty,\frac12]$ and equal $1$ on $[1,\infty)$; it is not the cut-off $\chi$ of the Toeplitz
operator. We write
$I_{\psi}(a):=\int_{0}^{\infty}e^{\sqrt{-1}t\psi}\kappa(t)\,a\,dt$, and
$I:=I_{\phi}$. We write $A\equiv B$ if $A-B$ is smooth near the diagonal.
The \emph{Reeb-independent class} $\mathcal C$ consists of the classical symbols
whose terms depend on $(z,w)$ only, and the class $\mathcal C'\supset\mathcal C$ of
those whose terms depend on $(z,w,\vartheta_{y})$, that is, do not depend on
$\vartheta_{x}$. For functions of $(z,w,\vartheta_{y})$, flatness means
$O(|z-w|^{N})$ for every $N$, locally uniformly in $\vartheta_{y}$. More generally, a function of
$(z,w,\vartheta_{x},\vartheta_{y})$ is \emph{flat at $z=w$} if it is $O(|z-w|^{N})$ for every $N$,
locally uniformly in $(\vartheta_{x},\vartheta_{y})$; this is stronger than flatness at $x=y$.

\paragraph{Inputs.}
\begin{enumerate}[label=\rm(H\arabic*)]
\item\label{hyp:H1} \cite[Thm.~1.9 with $q=0$, Thm.~5.1]{HM17JDG}: $\Pi\equiv I_{\phi_{-}}(s_{-})$,
      with $s_{-}$ classical of order $n$, and a phase $\phi_{-}$ satisfying ($\alpha$)--($\gamma$) of Theorem~\ref{thm:BdMS}.
      Here $\overline{\partial}_{b}$ has closed range \cite[Thm.~1.12]{HM17JDG}.
\item\label{hyp:H2} \cite[Rem.~3.6]{HM14}, which refers to \cite[Part~I, Chaps.~3, 7]{Hs10} and
      \cite{MenSj78}: a phase $\psi$ satisfying ($\alpha$)--($\gamma$) of Theorem~\ref{thm:BdMS}
      is equivalent to $\phi_{-}$ in the sense of Melin--Sj\"ostrand
      \cite[p.~172]{MeSj74}, that is, an oscillatory integral with phase $t\phi_{-}$
      and a classical symbol can be rewritten, modulo smoothing operators,
      as one with phase $t\psi$, and conversely. The reason is that the
      eikonal equation, with the first- and second-order terms prescribed
      by ($\alpha$) and ($\gamma$), determines the phase up to
      equivalence. Equivalence is characterized in \cite[Def.~4.3, Thm.~4.4]{HM17JDG}: a phase
      satisfying ($\alpha$) and the lower bound in ($\gamma$) is equivalent to $\phi_{-}$ if and only
      if it is a smooth multiple of $\phi_{-}$ modulo a function vanishing to infinite order on the
      diagonal; see also \cite[\S\S4, 5, 8]{HM17JDG}.
\end{enumerate}

We use \ref{hyp:H2} through Proposition~\ref{nf:prop:phase}.

\paragraph{Hypotheses.} Let $D$ be a BRT chart, let $dV$ be a smooth positive density on $D$ invariant
under the Reeb flow, and let $K$ be the distribution kernel, with respect to $dV$, of an operator on
$D$. We consider the following conditions.
\begin{enumerate}[label=\rm(K\arabic*)]
\item\label{K:rep} $K\equiv I(\sigma)$ near the diagonal, with $\sigma$ classical of order $n$ and
      $\sigma_{0}(x,x)\ne0$.
\item\label{K:inv} $K(e^{uT}x,e^{uT}y)=K(x,y)$ for small $u$.
\item\label{K:CR} $\overline{Z}_{l,x}K=0$ for $l=1,\dots,n$.
\item\label{K:sym} $K(x,y)=\overline{K(y,x)}$.
\end{enumerate}
For the Szeg\H o projection $\Pi$ of a compact Sasakian manifold and $dV=dV_{\eta}$: \ref{K:rep} holds
by \ref{hyp:H1} and Proposition~\ref{nf:prop:phase} below, the leading coefficient on the diagonal being
$1/(2\pi^{n+1})$ \cite[Thm.~4.5]{HM17JDG}, \cite[(1.1.4)]{CHH25}; \ref{K:inv} holds by
Lemma~\ref{lem:sa:commute}; \ref{K:CR} because $\Pi$ maps into CR functions; and \ref{K:sym} because
$\Pi$ is self-adjoint.

\subsection{Passage to the BRT phase}\label{nf:sec:phase}
\begin{proposition}[The BRT phase]\label{nf:prop:phase}
Let $\psi$ satisfy \textup{($\alpha$)--($\gamma$)} of Theorem~\ref{thm:BdMS} on a BRT chart, and let
$K\equiv I_{\psi}(\sigma)$ with $\sigma$ classical of order $n$. Then $K\equiv I(\sigma')$ for a classical
symbol $\sigma'$ of order $n$ with $\sigma'_{0}(x,x)=\sigma_{0}(x,x)$. More precisely, $\psi=g\phi+f$ with $g$
smooth, $g(x,x)=1$ and $f$ flat, and one can take $\sigma'_{j}=g^{-(n+1-j)}\sigma_{j}$.
\end{proposition}

\subsection{Reeb-independent normal form and uniqueness}\label{nf:sec:normal}
\begin{lemma}[Model integrals]\label{cf:lem:J}
For $k\in\mathbb Z$ and $\operatorname{Re}x>0$ let
$J_{k}(x)=\int_{0}^{\infty}e^{-tx}\kappa(t)t^{k}\,dt$. Then, with entire $E_{k}$
and the principal branch of $\log$,
\begin{equation}
 J_{k}(x)=k!\,x^{-k-1}+E_{k}(x)\quad(k\ge0),\qquad
 J_{-1-p}(x)=\tfrac{(-1)^{p+1}}{p!}\,x^{p}\log x+E_{-1-p}(x)\quad(p\ge0).
\end{equation}
\end{lemma}
\begin{proof}
For $k\ge0$, $J_{k}-k!x^{-k-1}=-\int_{0}^{\infty}e^{-tx}(1-\kappa)t^{k}dt$ is entire. For
$k=-1-p$, $\partial_{x}^{p+1}J_{-1-p}=(-1)^{p+1}J_{0}$ and
$\partial_{x}^{p+1}(x^{p}\log x)=p!\,x^{-1}$, so the difference of the two
sides has an entire $(p+1)$-st derivative. This is \cite[(3.1)]{HS20}.
\end{proof}

\begin{lemma}[Flat terms are negligible]\label{cf:lem:flat}
Call a smooth function $\psi$ near the diagonal a \emph{Szeg\H o phase} if $\psi(x,x)=0$,
$d_{x}\psi(x,x)=\eta(x)$ and $\operatorname{Im}\psi\ge c|z-w|^{2}$ for some $c>0$.
\textup{(a)} If $a$ is flat, then $I_{\psi}(a)$ is smooth near the diagonal,
for every Szeg\H o phase $\psi$. \textup{(b)} If $\psi,\psi'$ are Szeg\H o
phases with $\psi-\psi'$ flat, then $I_{\psi}(a)\equiv I_{\psi'}(a)$ for every
classical $a$.
\end{lemma}

\begin{proposition}[Uniqueness]\label{cf:prop:unique}
If $a\in\mathcal C'$ and $I(a)$ is smooth near the diagonal over $V$, then
every $a_{j}$ is flat over $V$.
\end{proposition}

\begin{proposition}[Reeb-independent normal form]\label{cf:prop:normal}
Let $K$ satisfy \ref{K:inv}. If $K\equiv I(\sigma)$ with $\sigma$ classical,
there is $\widehat{\sigma}\in\mathcal C$ with $K\equiv I(\widehat{\sigma})$. If $\sigma$
has no terms of negative order, neither has $\widehat{\sigma}$.
\end{proposition}

\subsection{Almost sesquiholomorphy and polarization}\label{cf:sec:CR}
From now on $K$ satisfies \ref{K:rep}--\ref{K:sym}, and $s\in\mathcal C$ is the symbol given by Proposition~\ref{cf:prop:normal}; by
Proposition~\ref{cf:prop:unique} it is unique modulo flat terms.

\begin{proposition}\label{cf:prop:CR}
For all $j$ and $l$, the functions $\partial_{\bar z_{l}}s_{j}$ and
$\partial_{w_{l}}s_{j}$ are flat.
\end{proposition}

\begin{corollary}[Polarization]\label{cf:cor:polar}
Each $s_{j}$ is determined modulo flat functions by $z\mapsto s_{j}(z,z)$; in
particular, if $s_{j}(z,z)=0$ on $V$, then $s_{j}$ is flat over $V$.
\end{corollary}

\subsection{Summary}\label{nf:sec:summary}
\begin{proposition}[Symbol normal form]\label{nf:prop:symbol}
Let $K$ satisfy \ref{K:rep}--\ref{K:sym}. Then $K\equiv I(s)$ with $s\in\mathcal C$ of order $n$, unique
modulo flat terms; $\partial_{\bar z}s_{j}$ and $\partial_{w}s_{j}$ are flat, and $s_{j}$ is determined modulo
flat functions by $z\mapsto s_{j}(z,z)$. This applies to the Szeg\H o projection $\Pi$ of a compact
Sasakian manifold, and the values $s_{j}(x,x)$ are then well defined functions on $X$.
\end{proposition}

\begin{remark}[Where the normal forms come from]\label{nf:rem:mechanism}
We include the proofs (Appendix~\ref{app:nf}) for two reasons. They are local, so they apply
to $(n,0)$-forms on complete manifolds (Section~\ref{tf:sec:setting}), where compactness is not
available. And they show, for readers from K\"ahler geometry, where the normal forms come from: the
phase is affine in the Reeb variable (P1), so that evaluation at its complex zero removes the dependence
on $\vartheta_{x}$ (Lemma~\ref{cf:lem:AA}); the invariance \ref{K:inv} removes that on $\vartheta_{y}$; the
Fourier transform in $\vartheta_{x}-\vartheta_{y}$ separates the orders (Proposition~\ref{cf:prop:unique});
and the Cauchy--Riemann equations \ref{K:CR} with the symmetry \ref{K:sym} give the almost
sesquiholomorphy (Proposition~\ref{cf:prop:CR}).
\end{remark}

\section{The coefficients: proof of Theorem~\ref{thm:introD}}\label{sec:coef}
This section carries out the additional step for functions. Proposition~\ref{prop:sa:exact} is the form
of the expansion of Section~\ref{ssec:sa:known} that we need. Its structure is in \cite[Thm.~1.3]{CHH25}:
when $T$ is a CR vector field, for a phase $\phi$ with $T_{x}\phi=1$ and $T_{y}\phi=-1$ the semi-classical
amplitude is $\chi(t)$ times a symbol independent of the Reeb variables. Proposition~\ref{prop:sa:exact}
identifies this symbol with the Reeb-independent symbol $s$ of $\Pi$ (Proposition~\ref{nf:prop:symbol}),
so that the coefficients are the values $s_{j}(x,x)$. The mechanism is elementary: in BRT coordinates the
flow of $T$ is the translation $\vartheta\mapsto\vartheta+u$, the BRT phase \eqref{cf:eq:BRTphase} is
$\vartheta_{x}-\vartheta_{y}+\mathcal B(z,w)$, and $s$ does not depend on $\vartheta$; hence
$\chi(-\sqrt{-1}T/k)$ acts on the oscillatory integral of $\Pi$ as the multiplier $\chi(\tau/k)$ in the
frequency $\tau$. The proof is an exact Fourier transform in the Reeb variable: no stationary phase is
needed, and the error terms come only from the localization. Lemma~\ref{lem:sa:transverse} identifies the
values $s_{j}(x,x)$, by comparison with a circle bundle over $\CP^{n}$.

\begin{proposition}[Exact identification]\label{prop:sa:exact}
In the coordinates of Lemma~\ref{lem:sa:BRT}, near the diagonal,
\begin{equation}\label{eq:sa:exact}
 \big(\chi(-\sqrt{-1}T/k)\Pi\big)(x,y)=\int_{0}^{\infty}e^{\sqrt{-1}\tau\phi(x,y)}\chi(\tau/k)\,s(z,w,\tau)\,d\tau
 +O(k^{-\infty}),
\end{equation}
and $\chi(-\sqrt{-1}T/k)\Pi$ is $O(k^{-\infty})$ away from the diagonal. In particular,
uniformly on $X$,
\begin{equation}\label{eq:sa:diag}
 \big(\chi(-\sqrt{-1}T/k)\Pi\big)(x,x)\sim\sum_{j\ge0}m_{n-j}(\chi)\,s_{j}(x,x)\,k^{n+1-j}.
\end{equation}
\end{proposition}
\begin{proof}
Let $\varrho\in \mathscr{C}^{\infty}_{c}((-\delta,\delta))$ equal $1$ near $0$, and split
\eqref{eq:sa:chiT} into $F_{\rm loc}=\int\varrho K_{k}\,e^{uT}\Pi\,du$ and $F_{\rm far}$.
Then $F_{\rm far}=G_{k}(-\sqrt{-1}T)\Pi$ with
$G_{k}(\lambda)=\chi(\lambda/k)-\frac1{2\pi}\int\widehat{\varrho}(\mu)\chi((\lambda-\mu)/k)\,d\mu$.
Since $\varrho\equiv1$ near $0$, all moments of $\widehat{\varrho}$ of positive order
vanish, and Taylor's formula in $\mu/k$ together with the decay of $\widehat{\varrho}$
gives $|G_{k}(\lambda)|\le C_{N}k^{-N}(1+|\lambda|)^{-N}$ for every $N$. On $H^{0}_{b}(X)$ the
operator $\Box_{b}-T^{2}=-T^{2}$ is elliptic, so powers of $1-\sqrt{-1}T$ control
all Sobolev norms there, and the kernel of $F_{\rm far}$ is $O(k^{-\infty})$ in
$\mathscr{C}^{\infty}$. No information on the returns of the flow of $T$ is needed.

For $(x,y)$ near the diagonal and $|u|<\delta$, $e^{uT}$ translates $\vartheta_{x}$ by
$u$ and preserves the chart. Write $\Pi=I(s)+R$ with $R$ smooth. The moments
of $\varrho K_{k}$ are $O(k^{-\infty})$, because those of $K_{k}$ vanish ($\chi=0$
near $0$) and $K_{k}=O(k^{-\infty})$ on the support of $1-\varrho$; so $R$
contributes $O(k^{-\infty})$, by Taylor's formula in $u$. Away from the
diagonal, let $(x,y)$ range over a compact set disjoint from it. Then
there is $c>0$ such that $e^{uT}x\ne y$ for $|u|<c$, so $u\mapsto\Pi(e^{uT}x,y)$ is
smooth on $(-c,c)$, uniformly, and this part contributes $O(k^{-\infty})$ by the
same moment argument. On $c/2\le|u|\le\delta$ the function $\Pi(e^{uT}x,y)$ may be
singular, at the points with $e^{uT}x=y$; but it is a distribution in
$u$ of finite order depending smoothly on $(x,y)$, because the wave front
set of $\Pi$ consists of multiples of $(\eta,-\eta)$ over the diagonal and
$\eta(T)=1$, and $K_{k}$ with all its derivatives is
$O(k^{-\infty})$ there. Hence $F_{\rm loc}$ is $O(k^{-\infty})$ away from the
diagonal. The contribution of $I(s)$ is
$\int e^{\sqrt{-1}\tau\phi}(\varrho K_{k})^{\vee}(\tau)\,s\,d\tau$ with
$(\varrho K_{k})^{\vee}(\tau)=\chi(\tau/k)+O(k^{-N}(1+\tau)^{-N})$, by the computation
above. This gives \eqref{eq:sa:exact}, and \eqref{eq:sa:diag} follows by
substituting $\tau=kt$.
\end{proof}

\begin{lemma}[Realization of a potential]\label{lem:sa:real}
Let $\varphi$ be smooth and strictly plurisubharmonic near $0\in\mathbb C^{n}$ with
$\varphi(z)=|z|^{2}+O(|z|^{3})$. Then there are $r>0$ and a Hermitian metric $h$
of positive curvature on $\mathcal O(1)\to\CP^{n}$ whose weight, in the affine chart
$\mathbb C^{n}\subset\CP^{n}$ and the standard frame, equals $\varphi$ on $\{|z|<r/2\}$.
\end{lemma}
\begin{proof}
In the affine chart the Fubini--Study metric has the weight
$\varphi_{\FS}=\log(1+|z|^{2})=|z|^{2}+O(|z|^{4})$, so $d:=\varphi-\varphi_{\FS}$ satisfies
$\partial^{\alpha}d=O(|z|^{3-|\alpha|})$ for $|\alpha|\le2$. Let $\kappa_{1}\in \mathscr{C}^{\infty}_{c}(\{|z|<1\})$
equal $1$ on $\{|z|\le\frac12\}$ and put $\psi_{r}(z)=\kappa_{1}(z/r)\,d(z)$, extended
by zero to $\CP^{n}$. On $\{|z|\le r\}$ the second derivatives of $\psi_{r}$
are bounded by $C(r^{-2}r^{3}+r^{-1}r^{2}+r)=O(r)$, while
$\sqrt{-1}\partial\overline{\partial}\varphi_{\FS}$ is bounded below by a positive multiple of the
Euclidean form on the unit ball. Hence $h=h_{\FS}e^{-\psi_{r}}$ has positive
curvature for $r$ small, and its weight $\varphi_{\FS}+\psi_{r}$ equals $\varphi$ on
$\{|z|<r/2\}$.
\end{proof}

\begin{lemma}[The universal coefficients for the transverse metric]\label{lem:sa:transverse}
Let $a_{j}(\cdot)$ be the universal polynomials in the curvature and its covariant derivatives with
$a^{L}_{j}=a_{j}(g)$ in \eqref{eq:intro:TYZ}. Then $s_{j}(x,x)=a_{j}(g^{T})(x)/(2\pi^{n+1})$ for every $j$.
\end{lemma}
\begin{proof}
Fix $x_{0}\in X$. By Lemma~\ref{lem:sa:chart} there is a chart $(z,\vartheta)\in U\times(-\varepsilon,\varepsilon)$
with $x_{0}=(0,0)$ and potential $\varphi(z)=|z|^{2}+O(|z|^{3})$. Let $h$ and $r$ be as in
Lemma~\ref{lem:sa:real}, and let $X'$ be the unit circle bundle of $(\mathcal O(-1),h^{*})$ over
$\CP^{n}$, with its connection form $\eta'$; it is a regular compact Sasakian manifold. Over
$\{|z|<r/2\}$ and for $|\vartheta|<\varepsilon$, the coordinates $(z,\vartheta)$ on $X'$ form a chart as in
Lemma~\ref{lem:sa:BRT} with the same potential $\varphi$, so the identity map between the two charts is
a CR diffeomorphism which preserves the contact forms and the volume forms. By \eqref{eq:sa:diag}, the
values $s_{j}(x,x)$ and $s'_{j}(x,x)$ are the coefficients $a_{j}$ of \cite[Thm.~1.5]{CHH25} for
$(X,\eta,dV_{\eta})$ and $(X',\eta',dV_{\eta'})$, so $s_{j}(x_{0},x_{0})=s'_{j}(x_{0},x_{0})$ by
\cite[Rem.~1.6]{CHH25}. On the circle bundle $X'$, where $dV_{\eta'}=d\theta\wedge\pi^{*}dV_{g}$,
\cite[Ex.~1.7]{CHH25}, that is \cite[Thm.~4.2.3]{HHMS23}, gives $s'_{j}(x,x)=b_{j}(\pi(x))/(2\pi)$, where
$B^{L}_{k}\sim\sum_{j}b_{j}k^{n-j}$; by \eqref{eq:intro:TYZ}, $b_{j}=a^{L}_{j}\pi^{-n}$, so
$s'_{j}(x_{0},x_{0})=a^{L}_{j}(0)/(2\pi^{n+1})$. Finally $a^{L}_{j}=a_{j}(g)$, and on $\{|z|<r/2\}$ the metric
$g$ has the K\"ahler form $\frac{\sqrt{-1}}{2}\partial\overline{\partial}\varphi=\frac12d\eta$, which is the
transverse K\"ahler form of $X$ in the chart. Hence $a_{j}(g)(0)=a_{j}(g^{T})(x_{0})$.
\end{proof}

\begin{remark}[Where the line bundle enters]\label{rem:sa:model}
The statement of Lemma~\ref{lem:sa:transverse} involves no line bundle, and the proof uses one only on the
comparison model $X'$, never on $X$ or on a quotient of $X$. The passage from the Toeplitz coefficients to
the values $s_{j}(x,x)$, Proposition~\ref{prop:sa:exact}, needs only BRT coordinates and the normal form:
it is an exact Fourier transform in the Reeb variable.
For a circle bundle the identification of $s_{j}(x,x)$ is also a computation in BRT coordinates: by
\eqref{eq:fourier} and the Reeb-independence of $s$, $S_{k}(x,x)=s(z,z,k)+O(k^{-\infty})$, as in the proof of
Proposition~\ref{tf:prop:diag}, while $S_{k}(x,x)=B^{L}_{k}/(2\pi)$ by \eqref{eq:SkBk}. For a general $X$, BRT
coordinates show only that $s_{j}(x,x)$ is a universal expression in the jets of $\varphi$ at $x$
\cite[Rem.~3.9]{CHH25}, and three facts identify it with $a_{j}(g^{T})(x)/(2\pi^{n+1})$. First, the
coefficients are local \cite[Rem.~1.6]{CHH25}. Second, whether $X$ is regular, quasi-regular or irregular is
a property of the closures of the Reeb orbits, which a germ does not see: in the flow-box coordinates of
Lemma~\ref{lem:sa:BRT} the Sasakian structure is determined by $\varphi$, so every germ of a Sasakian
manifold is the germ of a circle bundle over a ball. Third, Lemma~\ref{lem:sa:real} makes this circle bundle
compact and regular, as required by \cite{CHH25} and \cite{HHMS23}, which concern compact manifolds. Some
comparison with line bundles cannot be avoided, since the polynomials $a_{j}(\cdot)$ are defined through
them; any compact model realizing the germ would do, and $\CP^{n}$ is the simplest.
Theorem~\ref{tf:thm:id} uses the same argument, with the twisted coefficients.
\end{remark}

\begin{proof}[Proof of Theorem~\ref{thm:introD}]
By Proposition~\ref{nf:prop:symbol}, $\Pi\equiv I(s)$ with $s$ Reeb-independent and unique modulo flat
terms, and $s_{j}$ is determined modulo flat functions by its diagonal values.
Proposition~\ref{prop:sa:exact} gives the expansion, with coefficients $m_{n-j}(\chi)s_{j}(x,x)$, and
Lemma~\ref{lem:sa:transverse} identifies them.
\end{proof}

\section{The logarithmic term: proof of Theorem~\ref{thm:introB}}\label{sec:log}
\begin{definition}\label{cf:def:log}
The logarithmic term of $\Pi$ \emph{vanishes to infinite order over} an open set $W\subset X$ if, near
each diagonal point over $W$,
\[
 \Pi\equiv u\,(-\sqrt{-1}\psi)^{-(n+1)}+v\log(-\sqrt{-1}\psi),
\]
where $\psi$ is a phase satisfying \textup{($\alpha$)--($\gamma$)} of Theorem~\ref{thm:BdMS}, hence equivalent to
the BRT phase by \ref{hyp:H2}, and $u,v$ are smooth with $v$ flat.
\end{definition}
A flat $v$ can be absorbed into $u$, since $v\log(-\sqrt{-1}\psi)$ is then smooth. By
Lemma~\ref{cf:lem:J}, for $s\in\mathcal C$,
\begin{equation}\label{eq:log:G}
\begin{gathered}
 I(s)\equiv\sum_{j\le n}(n-j)!\,s_{j}\,(-\sqrt{-1}\phi)^{j-n-1}+G\log(-\sqrt{-1}\phi),\\
 G=\sum_{p\ge0}\frac{(-1)^{p+1}}{p!}\,s_{n+1+p}\,(-\sqrt{-1}\phi)^{p},
\end{gathered}
\end{equation}
where the sum defining $G$ is understood as a Borel sum.

\begin{theorem}[The logarithmic term]\label{thm:logS}
Let $(X,\eta,T)$ be a compact Sasakian manifold, let $W\subset X$ be open, and let $s$ be the symbol of
Proposition~\ref{nf:prop:symbol} for $\Pi$. The following are equivalent:
\begin{enumerate}[label=\rm(\roman*)]
\item the logarithmic term of $\Pi$ vanishes to infinite order over $W$;
\item $s_{n+1+p}$ is flat over $W$ for every $p\ge0$;
\item $a_{n+1+p}(g^{T})=0$ on $W$ for every $p\ge0$;
\item for some $\chi\ge0$, $\chi\ne0$, the coefficients of $k^{-p}$, $p\ge0$, in \eqref{eq:sa:diag} vanish on $W$.
\end{enumerate}
\end{theorem}
\begin{proof}
(ii)$\Rightarrow$(i). By Lemma~\ref{cf:lem:flat}, $\Pi\equiv I(\sum_{j\le n}s_{j}t^{n-j})$, and by
\eqref{eq:log:G} this is a decomposition with $\psi=\phi$ and $v=0$.

(i)$\Rightarrow$(ii). Let $\Pi\equiv u(-\sqrt{-1}\psi)^{-(n+1)}+v\log(-\sqrt{-1}\psi)$ near a diagonal point
over $W$, with $v$ flat. Then $v\log(-\sqrt{-1}\psi)$ is smooth, and by Lemma~\ref{cf:lem:J} and
Proposition~\ref{nf:prop:phase}, with $\psi=g\phi+f$ as there,
$\Pi\equiv I_{\psi}(\frac{u}{n!}t^{n})\equiv I(\sigma)$ with $\sigma=g^{-(n+1)}\frac{u}{n!}t^{n}$. This symbol has no
terms of negative order, so Proposition~\ref{cf:prop:normal} gives $\widehat\sigma\in\mathcal C$ without terms of
negative order and with $\Pi\equiv I(\widehat\sigma)$. Hence $I(s-\widehat\sigma)$ is smooth near the diagonal,
and Proposition~\ref{cf:prop:unique} shows that $s_{j}-\widehat\sigma_{j}$ is flat for every $j$. For $j\ge n+1$,
$\widehat\sigma_{j}=0$, so $s_{j}$ is flat.

(ii)$\Leftrightarrow$(iii). By Theorem~\ref{thm:introD}, $s_{n+1+p}(x,x)=a_{n+1+p}(g^{T})(x)/(2\pi^{n+1})$;
so (ii) implies (iii), and Corollary~\ref{cf:cor:polar} gives the converse.

(iii)$\Leftrightarrow$(iv). By \eqref{eq:sa:diag}, the coefficient of $k^{-p}$ is
$m_{-1-p}(\chi)s_{n+1+p}(x,x)$, and $m_{-1-p}(\chi)>0$ for such $\chi$.
\end{proof}

\begin{remark}[Comparison with \cite{LT}]\label{rem:logLT}
(a) On the interior diagonal of the disc bundle of a circle bundle, $-\sqrt{-1}\phi=-\log t$ with
$t=|v|^{2}_{h^{*}}$, and the $p$-th term of $G$ in \eqref{eq:log:G} is the logarithmic part of
$\operatorname{Li}_{1+p}(t)$ (Lemma~\ref{lem:sums}).

(b) For circle bundles, the phase $\psi_{\rm LT}$ of \cite[(2.1)]{LT} equals
$-\sqrt{-1}(e^{\sqrt{-1}\phi}-1)$ on $X\times X$ modulo functions flat at $z=w$, since the almost analytic
extension of the metric used in \cite[\S2]{LT} and $e^{-\varphi(z,\overline w)}$ differ by such functions;
as $-\sqrt{-1}(e^{\sqrt{-1}\phi}-1)=\phi\,g$ with $g$ smooth and $g(x,x)=1$, it satisfies \textup{($\alpha$)--($\gamma$)} and is equivalent to $\phi$.
Since $\log\psi=\log(-\sqrt{-1}\psi)+\frac{\pi}{2}\sqrt{-1}$ for $\operatorname{Im}\psi\ge0$, the vanishing of
the coefficient $v$ in \cite[(2.4)]{LT}, which Lu and Tian require on $\overline D\times\overline D$, implies (i)
with $W=X$. Conversely, under (iii) the proof of (ii)$\Rightarrow$(i), with $\phi$ written as
$\psi_{\rm LT}/g$, gives a decomposition of the form \cite[(2.4)]{LT} on $X\times X$ near the diagonal with
$v=0$. Definition~\ref{cf:def:log} is local and concerns only $X\times X$.

(c) For circle bundles, by Proposition~\ref{prop:dictionaryS}, (iii) is also equivalent to the flatness at
$X$ of the logarithmic coefficient $\Psi_{S}$ of the Szeg\H o kernel on the diagonal of $D$, also over
open subsets of $M$.
\end{remark}

\begin{example}[The sphere, off the diagonal]\label{cf:ex:sphere}
Let $M=\mathbb{CP}^{n}$ and $L=\mathcal O(1)$, so that $X=S^{2n+1}$. In the
affine chart, with $\varphi=\log(1+|z|^{2})$ and the lift
$x=e^{\sqrt{-1}\vartheta}(1,z)/(1+|z|^{2})^{1/2}$, one has
$\varphi(z,\overline{w})=\log(1+z\cdot\overline{w})$, hence
\begin{equation}
 \langle x,y\rangle_{\mathbb C^{n+1}}
 =e^{\sqrt{-1}(\vartheta_{x}-\vartheta_{y})}\,
 \frac{1+z\cdot\overline{w}}{(1+|z|^{2})^{1/2}(1+|w|^{2})^{1/2}}
 =e^{\sqrt{-1}\phi(x,y)}
\end{equation}
exactly. So the Szeg\H o kernel of the sphere, $\frac{n!}{2\pi^{n+1}}(1-\langle x,y\rangle)^{-(n+1)}$, reads
$\Pi=\frac{n!}{2\pi^{n+1}}\big(1-e^{\sqrt{-1}\phi}\big)^{-(n+1)}
=\frac{1}{2\pi^{n+1}}\sum_{k\ge0}\prod_{i=1}^{n}(k+i)\,e^{\sqrt{-1}k\phi}$,
and the Poisson summation formula \cite[Thm.~7.2.1]{Hoe03}, applied to
$t\mapsto\rho(t)\prod_{i}(t+i)\,e^{\sqrt{-1}t\phi}$ with $\rho\in \mathscr{C}^{\infty}(\mathbb R)$ equal to $1$
on $[0,\infty)$ and to $0$ on $(-\infty,-\frac12]$, gives near the diagonal
\begin{equation}
 \Pi\equiv I(s),\qquad s(t)=\frac{1}{2\pi^{n+1}}\prod_{i=1}^{n}(t+i):
\end{equation}
the terms with $\ell\ne0$ in the formula are smooth there, since
$|\operatorname{Re}\phi|<\pi$, and their sum converges since $n\ge1$.
Thus the phase is exactly $\phi$; the symbol is Reeb-independent with constant,
hence sesquiholomorphic, coefficients; and $s_{n+1+p}=0$ for all $p$. This agrees with Proposition~\ref{nf:prop:symbol} and Theorem~\ref{thm:logS}, since $a^{L}_{j}=0$ for $j\ge n+1$ on
$\mathbb{CP}^{n}$.
\end{example}

\section{Circle bundles and the disc bundle}\label{sec:circle}\label{app:log}
\emph{Setting.} In this section $M$ is compact, $(L,h)$ is positive and $X$ is the circle
bundle of $L^{*}$, so that the Reeb flow is the circle action. Only the identification of the Fourier components of $S$ with $H^{0}(M,L^{k})$ and the semi-classical expansion of $B^{L}_{k}$ are used; the Szeg\H o parametrix plays no
role.

By \eqref{eq:SzegoD} below, the Bergman kernel of $D$
\cite[Thm.~2.2]{EXX} and the Szeg\H o kernel of $D$ are, on the diagonal and
up to explicit factors, power series in $t=|v|^{2}_{h^{*}}$ with coefficients
$kB_{k}$ and $B^{L}_{k}$. The relation
between their logarithmic terms at $X$ and the coefficients of the
semi-classical expansion \eqref{eq:intro:TYZ} is therefore a question about
the singularity of a power series at $t=1$. This section derives the resulting formulas for the logarithmic terms.

\subsection{The two kernels}\label{sec:twokernels}

This section uses the kernels attached to $X$ and
to $D$, which we fix first.

Let $H^{0}_{b}(X)\subset L^{2}(X)$ be the Hardy space and $S$ the Szeg\H o
projector, both taken with respect to $dV_{X}=d\theta\wedge\pi^{*}dV_{g}$,
which is the measure $d\mu$ of \cite[\S2]{LT}. With this measure, Lu and Tian
write $\Pi_{m}=2\pi S_{m}$ \cite[(2.2)]{LT}, so that their $\Pi_{m}(x,x)$ is our
$B^{L}_{m}$; see \eqref{eq:SkBk}.

The $S^{1}$-action is the Reeb flow of $\eta$, so it commutes with $S$ by
Lemma~\ref{lem:sa:commute}. It gives
$H^{0}_{b}(X)=\widehat\bigoplus_{k\ge0}H^{0}_{b,k}(X)$ with
$H^{0}_{b,k}(X)\cong H^{0}(M,L^{k})$, and $S=\sum_{k}S_{k}$ with
\begin{equation}\label{eq:fourier}
 S_{k}(x,y)=\frac{1}{2\pi}\int_{-\pi}^{\pi}e^{-\sqrt{-1}k\theta}\,
 S\big(e^{\sqrt{-1}\theta}\!\cdot\! x,\,y\big)\,d\theta .
\end{equation}
Writing $B^{L}_{k}$ for the Bergman kernel function of $(L^{k},h^{k})$ with
respect to $dV_{g}$, fiber integration gives
\[
 \int_{X}|\widehat{s}|^{2}dV_{X}=2\pi\int_{M}|s|^{2}_{h^{k}}dV_{g},
\]
hence
\begin{equation}\label{eq:SkBk}
 S_{k}(x,x)=\tfrac{1}{2\pi}B^{L}_{k}(\pi(x)) .
\end{equation}
For $v\in D$ with $t=|v|^{2}_{h^{*}}$ one has
$|\widehat{s}(v)|^{2}=|s|^{2}_{h^{k}}t^{k}$, so the Szeg\H o kernel of $D$ is
\begin{equation}\label{eq:SzegoD}
 S_{D}(v,v)=\frac{1}{2\pi}\sum_{k\ge0}B^{L}_{k}(z)\,t^{k},
 \qquad\text{versus}\qquad
 \mathcal K=e^{\varphi}\det g\cdot\frac1\pi\sum_{k\ge1}kB_{k}(z)\,t^{k-1}.
\end{equation}
In the left member, the term with index $j$ of
$B^{L}_{k}\sim(k/\pi)^{n}\sum_{j}a^{L}_{j}k^{-j}$ contributes the exponent $k^{n-j}$.
By Lemma~\ref{lem:sums}, the logarithmic term of $S_{D}$ therefore comes
from the indices $j\ge n+1$, and that of $\mathcal K$, whose coefficients carry
the extra factor $k$, from the indices $j\ge n+2$. The twist by $K_{M}$
changes the coefficients, from $a^{L}_{j}$ to $a^{E}_{j}$, but not this
threshold.

\subsection{Summation}

\begin{lemma}\label{lem:sums}
Let $t\in(0,1)$ and $\rho=1-t$.
\begin{enumerate}[label=\rm(\alph*)]
\item For $m\ge0$, $\displaystyle\sum_{k\ge1}k^{m}t^{k}
      =\frac{A_{m}(t)}{\rho^{\,m+1}}$, where $A_{m}$ is a polynomial with
      $A_{m}(1)=m!$ ($A_{0}(t)=t$, and for $m\ge1$, $A_{m}(t)$ is $t$ times
      the Eulerian polynomial);
      in particular the sum is a rational function of $t$ with a pole of
      order $m+1$ at $t=1$ and no logarithmic term, and
      $\sum_{k\ge1}k^{m}t^{k}=m!\,\rho^{-(m+1)}\big(1+O(\rho)\big)$.
\item For $p\ge0$, with $\mu:=\log t$,
      \begin{equation}
       \sum_{k\ge1}k^{-1-p}t^{k}=\operatorname{Li}_{1+p}(t)
       =\frac{\mu^{p}}{p!}\big[H_{p}-\log(-\mu)\big]
       +\sum_{i\ne p}\frac{\zeta(1+p-i)}{i!}\,\mu^{i},
      \end{equation}
      valid for $|\mu|<2\pi$, where $H_{p}=\sum_{r\le p}1/r$. Since
      $\mu=-\rho+O(\rho^{2})$ and $\log(-\mu)=\log\rho+O(\rho)$, the
      logarithmic part of $\operatorname{Li}_{1+p}(t)$ is
      $-\frac{(-1)^{p}}{p!}\rho^{p}\big(1+O(\rho)\big)\log\rho$.
\end{enumerate}
\end{lemma}

Thus, in the expansion $kB_{k}\sim\pi^{-n}\sum_{j\ge0}a^{E}_{j}k^{\,n+1-j}$ of
the coefficients of
\begin{equation}
 \mathcal G(z,t):=\sum_{k\ge1}kB_{k}(z)\,t^{k-1}
 =\pi\,\mathcal K\,e^{-\varphi}(\det g)^{-1}
\end{equation}
the index $j$ contributes a pole of order $n+2-j$ and no logarithmic term
when $0\le j\le n+1$, and a term of type $\rho^{\,j-n-2}\log\rho$ when $j\ge n+2$.
Hence the logarithmic term of $\mathcal K$ begins with $j=n+2$.

\begin{remark}[The leading coefficient]\label{rem:leadingC}
The leading term of the Bergman kernel $\mathcal K$ of $D$ at $X$, namely
\begin{equation}
 \mathcal K=\frac{(n+1)!}{\pi^{n+1}}e^{\varphi}\det g\,\rho^{-(n+2)}(1+O(\rho)),
 \qquad\text{equivalently}\qquad \Phi|_{X}=\frac{(n+1)!}{2^{n+1}\pi^{n+1}},
\end{equation}
is
Lemma~\ref{lem:sums}(a) with $j=0$ and $m=n+1$: only
$B_{k}=(k/\pi)^{n}(1+O(k^{-1}))$ is needed, and the remainder contributes
$O(\rho^{-(n+1)})$.  For the ball,
$K_{\BB^{n+1}}=\frac{(n+1)!}{\pi^{n+1}}(1-|z|^{2})^{-(n+2)}$, in agreement with
this formula.
\end{remark}

\subsection{Logarithmic coefficients and the theorem of Lu--Tian}

The following proposition is \cite[Thm.~1.9]{EXX}, with the same
normalization and the same diagonal coefficients; its final equivalence is
\cite[Thm.~1.11, Cor.~1.12]{EXX}, and for compact $M$ the implication
``$\Rightarrow$'' is \cite[Thm.~2.4]{LT}. We include the summation proof,
which is short. The decomposition \eqref{eq:feff} is discussed in
Remark~\ref{rem:apriori}.

\begin{proposition}[Logarithmic coefficients of the Bergman kernel of $D$]\label{prop:dictionary}
Suppose that $\mathcal K$ admits, near $X$, the Fefferman-type
decomposition
\begin{equation}\label{eq:feff}
 \mathcal K=\Big(\Phi\,\rho^{-(n+2)}+\Psi\log\rho\Big)\,2^{n+1}e^{\varphi}\det g,
 \qquad \Phi,\Psi\in \mathscr{C}^{\infty}(\overline D)
\end{equation}
and write
$\Phi=\sum_{i}\phi_{i}\rho^{i}$, $\Psi=\sum_{p}\psi_{p}\rho^{p}$ for the
Taylor expansions at $X$; since $\mathcal K$ is invariant under the circle
action, $\phi_{i}$ and $\psi_{p}$ are functions on $M$. Then
$\phi_{0},\dots,\phi_{n+1}$ and all $\psi_{p}$ are linear combinations of the
$a^{E}_{j}$, with coefficients depending only on $n$ and the indices, and
\begin{equation}
 \psi_{p}=-\frac{(-1)^{p}}{2^{n+1}\pi^{n+1}\,p!}\;a^{E}_{n+2+p}
 \;+\;\Big(\text{a linear combination of }a^{E}_{n+2+p'},\ p'<p\Big),
\end{equation}
in particular
\begin{equation}
 \Psi\big|_{X}=\psi_{0}=-\frac{a^{E}_{n+2}}{2^{n+1}\pi^{n+1}}.
\end{equation}
Consequently $\Psi$ vanishes to infinite order at $X$ if and only if
$a^{E}_{n+2+p}=0$ for all $p\ge0$.
\end{proposition}

\begin{proof} This is the case of compact $M$ of \cite[Thm.~1.9]{EXX}; we include the short proof because it
parallels that of Proposition~\ref{prop:dictionaryS}.
Insert $kB_{k}=\pi^{-n}\big(\sum_{j\le N}a^{E}_{j}k^{\,n+1-j}+O(k^{\,n-N})\big)$
into \cite[Thm.~2.2]{EXX}. By Lemma~\ref{lem:sums} the explicit part
produces the poles $\rho^{-(n+2-j)}$ for $j\le n+1$ and the terms
$-\frac{(-1)^{p}}{p!}\rho^{p}(1+O(\rho))\log\rho$ for $j=n+2+p$, while the
remainder $\mathcal R_{N}(t):=\sum_{k}\epsilon_{k}t^{k-1}$, with
$\epsilon_{k}=O(k^{\,n-N})$, is of class $C^{N-n-2}$ up to $t=1$: its $M$-th
$t$-derivative has coefficients $O(k^{\,M+n-N})$, which are summable as soon
as $M\le N-n-2$. Since $\rho^{p}\log\rho$ is not of class $\mathscr{C}^{p}$ while
$\mathcal R_{N}$ is of class $C^{N-n-2}$, and since the polar part
$\Phi\rho^{-(n+2)}$ of \eqref{eq:feff} carries no logarithm,
$\mathcal R_{N}$ cannot contribute to $\psi_{p}$ for $p\le N-n-2$.
Comparing with \eqref{eq:feff} therefore determines $\psi_{p}$ for
$p\le N-n-2$, and letting $N\to\infty$ determines $\phi_{i}$ ($i\le n+1$)
and all $\psi_{p}$. The
factor $1/t=1+O(\rho)$ and the $O(\rho)$ inside
Lemma~\ref{lem:sums}(b) produce only lower-triangular corrections; the
diagonal coefficients $-(-1)^{p}/(2^{n+1}\pi^{n+1}p!)$ are non-zero, so the
correspondence is invertible.
\end{proof}

\begin{remark}[The role of \eqref{eq:feff}]\label{rem:apriori}
Both directions of Proposition~\ref{prop:dictionary} are Abelian: they use
only the summation of the expansion term by term, and no Tauberian
argument. Granted
the semi-classical expansion \eqref{eq:intro:TYZ} of $B_{k}$ to all orders
\cite{Cat99,Zel98,Lu00,MM07} and the decomposition \eqref{eq:feff}, the proof computes
the $\psi_{p}$ as linear functions of the $a^{E}_{j}$ by summing the expansion
term by term (Lemma~\ref{lem:sums}). The remainder is $C^{N-n-2}$ up to the
boundary, hence produces no $\rho^{p}\log\rho$ with $p\le N-n-2$. The
triangular system has non-zero diagonal entries, so both implications of
\begin{equation}
 a^{E}_{n+2+p}=0\ \forall p\ \Longleftrightarrow\ \Psi\ \text{flat at }X
\end{equation}
follow. The content of \eqref{eq:feff} is the existence of a decomposition
with $\Phi,\Psi$ smooth up to the boundary. For the diagonal kernel this also
follows from the summation, with the remainder estimates in $\mathscr{C}^{m}(M)$ for
every $m$ and Borel's lemma. Fefferman's
theorem \cite{Fe74} adds the identification of $\Psi$ with the CR-invariant logarithmic
term of a local model; for $D$ it is transported from a smoothly bounded
strongly pseudoconvex domain in $\mathbb C^{n+1}$ by the localization
\cite[Prop.~3.1]{HL}. The passage from the
global disc bundle to the germ at a boundary point, where Fefferman's
theorem does not apply, is treated in \cite[Thm.~1.2]{EXX}.
\end{remark}

The same computation, applied to the Szeg\H o kernel instead of the Bergman
kernel of $D$, gives the theorem of Lu and Tian. We record it because the proof below uses only Lemma~\ref{lem:sums} and the a priori structure.

For the Szeg\H o kernel of $D$, the decomposition \eqref{eq:feffS} below exists by
\eqref{eq:SzegoD} and the summation, as for \eqref{eq:feff} in Remark~\ref{rem:apriori}; for
smoothly bounded strongly pseudoconvex domains in $\mathbb C^{n+1}$ it goes back to
\cite[Thm.~1.5, Cor.~1.7]{BouSj75}.

\begin{proposition}[Logarithmic coefficients of the Szeg\H o kernel]\label{prop:dictionaryS}
Suppose that the Szeg\H o kernel of $D$ admits, near $X$, the decomposition
\begin{equation}\label{eq:feffS}
 S_{D}=\Phi_{S}\,\rho^{-(n+1)}+\Psi_{S}\log\rho,
 \qquad \Phi_{S},\Psi_{S}\in \mathscr{C}^{\infty}(\overline D)
\end{equation}
and write $\Psi_{S}=\sum_{p}\psi^{S}_{p}\rho^{p}$; as in
Proposition~\ref{prop:dictionary}, the $\psi^{S}_{p}$ are functions on $M$. Then
\begin{equation}
 \psi^{S}_{p}
 =-\frac{(-1)^{p}}{2\pi^{\,n+1}\,p!}\;a^{L}_{n+1+p}
 \;+\;\Big(\text{a linear combination of }a^{L}_{n+1+p'},\ p'<p\Big),
\end{equation}
with coefficients depending only on $n$ and $p$; in particular
$\Psi_{S}|_{X}=\psi^{S}_{0}=-a^{L}_{n+1}/(2\pi^{\,n+1})$.
\end{proposition}

\begin{proof}
By \eqref{eq:SzegoD}, $S_{D}(v,v)=\frac{1}{2\pi}\sum_{k\ge0}B^{L}_{k}(z)t^{k}$
with $B^{L}_{k}\sim\pi^{-n}\sum_{j}a^{L}_{j}k^{\,n-j}$. The argument is
verbatim that of Proposition~\ref{prop:dictionary}, with the exponent
$n+1-j$ there replaced by $n-j$ here: by Lemma~\ref{lem:sums}(a) the indices
$j\le n$ give poles $\rho^{-(n+1-j)}$ and no logarithmic term, and by
Lemma~\ref{lem:sums}(b) the index $j=n+1+p$ gives
$-\frac{(-1)^{p}}{p!}\rho^{p}(1+O(\rho))\log\rho$; the remainder after $N$
terms is of class $C^{N-n-1}$ up to $t=1$ and so contributes no
$\rho^{p}\log\rho$ with $p\le N-n-1$. Multiplying by $\frac{1}{2\pi}\pi^{-n}$
gives the displayed diagonal coefficient.
\end{proof}

We call $X$ \emph{Bergman logarithmically flat} if $\Psi$ in \eqref{eq:feff}
vanishes to infinite order at $X$, and \emph{Szeg\H o logarithmically flat} if
$\Psi_{S}$ in \eqref{eq:feffS} does.

The implication ``$\Rightarrow$'' in the following theorem is the counterpart, on the
diagonal of $D$, of \cite[Thm.~2.3]{LT}. Lu and Tian state the hypothesis off the diagonal, as the
vanishing of the logarithmic coefficient in the decomposition \cite[(2.4)]{LT} of the Szeg\H o kernel
on $\overline D\times\overline D$; restricted to the diagonal of $D$, where their phase equals
$\sqrt{-1}\rho$, this decomposition takes the form \eqref{eq:feffS} with $\Psi_{S}=0$. The
hypothesis below is the formulation on the diagonal used in \cite{ALZ13,LUZ16}. For expansions
that are finite in the exact sense the converse is \cite[Thm.~1.1]{LMZ19}; see also \cite{ALZ13}.
In general the converse is recorded in \cite{LMZ19} as an open conjecture, attributed to Lu; for
the Bergman kernel, Ebenfelt, Xiao and Xu describe the converse as a question posed by Lu and
Tian \cite[p.~6]{EXX} and prove it \cite[Thm.~1.11, Cor.~1.12]{EXX}. Below, both directions
on the diagonal of $D$ are obtained by summing \eqref{eq:intro:TYZ} along the fibers. On $X\times X$
near the diagonal, the converse follows in the same way from the off-diagonal expansion of the
Bergman kernels of $L^{k}$ \cite[Thm.~3.1, (3.5)]{BBS08}; Theorem~\ref{thm:logS} gives both
directions there in local form, in the sense of Definition~\ref{cf:def:log}, through the structure of
the Szeg\H o symbol, which does not involve the line bundle.

\begin{theorem}[The Lu--Tian correspondence, diagonal form]\label{thm:LuTian}
Let $M$ be compact, $L$ positive with metric
$h$, and $D$ the unit disc bundle of $(L^{*},h^{*})$. Then the coefficient
$\Psi_{S}$ in \eqref{eq:feffS} vanishes to infinite order at $X$ if and only
if
\begin{equation}
 a^{L}_{j}=0\qquad\text{for every }j\ge n+1 .
\end{equation}
\end{theorem}

\begin{proof}
The correspondence of Proposition~\ref{prop:dictionaryS} is triangular, with
diagonal entries
\begin{equation}
 -\frac{(-1)^{p}}{2\pi^{\,n+1}\,p!}\ne0 .
\end{equation}
Hence it is invertible, and all $\psi^{S}_{p}$ vanish if and only if all
$a^{L}_{n+1+p}$ do.
\end{proof}

\begin{remark}[The conjecture of Lu and Tian]\label{rem:LTcost}
Lu and Tian conjecture \cite[\S1]{LT} that on $\CP^{n}$, with the hyperplane
bundle $H$, a K\"ahler form $\omega\in[\omega_{\FS}]$ whose Szeg\H o kernel has
vanishing logarithmic term is $F^{*}\omega_{\FS}$ for an automorphism $F$; they
prove it for
$n=1$ \cite[Thm.~1.1]{LT} and near $\omega_{\FS}$ in all dimensions
\cite[Thm.~1.2]{LT}; local rigidity near the invariant metrics of integral K\"ahler--Einstein
polarizations of products of irreducible compact Hermitian symmetric spaces
is proved in \cite{LPZ26}. Corollary~\ref{cf:cor:rigid} gives a version in which, for real
analytic $\omega$, the hypothesis is imposed only over an open set. The
conjecture concerns $\CP^{n}$ only. For a flat torus, and for a compact
Hermitian symmetric space with a homogeneous line bundle \cite{EZ10}, the
argument of Remark~\ref{rem:models} gives $a^{L}_{j}=0$ for $j\ge n+1$, so
$X$ is Szeg\H o logarithmically flat by Proposition~\ref{prop:dictionaryS}.
\end{remark}

\begin{remark}[The models]\label{rem:models}
For each of the following models, $B_{k}$ is a polynomial
in $k$ of degree $n$ with leading coefficient $\pi^{-n}$:
\begin{equation}
 \CP^{r}(\mathcal O(q)):\ \frac{\prod_{j=1}^{r}(qk-j)}{q^{r}\pi^{r}},
 \qquad
 \BB^{l}(q):\ \frac{\prod_{j=1}^{l}(qk+j)}{q^{l}\pi^{l}},
 \qquad
 \text{flat}:\ \Big(\frac{k}{\pi}\Big)^{n},
\end{equation}
and for a product of models $B_{k}$ is the product of the factors.
Consequently $a_{0}=1$ and $a_{j}=0$ for all $j>n$, and $X$ is Bergman
logarithmically flat.

For $\CP^{r}$, the only compact model, the group of automorphisms of
$(M,g;L,h)$ acts transitively, so $B_{k}$ is constant and equals
$h^{0}(M,E_{k})/\Vol_{g}(M)$;
Kodaira vanishing and Riemann--Roch make $h^{0}(M,E_{k})$ a polynomial in
$k$, and the displayed formulas are elementary (for
$\CP^{r}$, $h^{0}(\mathcal O(qk-r-1))=\binom{qk-1}{r}$, which vanishes
together with the displayed product when $1\le qk\le r$). The non-compact
models $\BB^{l}(q)$ and $\mathbb C^{n}$ are homogeneous, and their $B_{k}$
are the weighted Bergman and Fock kernels, which give the displayed
products. If $B_{k}=\pi^{-n}P(k)$ with $\deg P=n$ then
$\sum_{j}a^{E}_{j}k^{-j}=P(k)/k^{n}$ is a polynomial in $k^{-1}$ of degree $n$.
Finally, $\sum_{k}kB_{k}t^{k-1}$ is then a rational function of $t$ by
Lemma~\ref{lem:sums}(a), so $\Psi=0$ in \eqref{eq:feff}.

\emph{Compact quotients.} For a compact locally homogeneous $(M,g;L,h)$
modeled on these spaces (a compact ball quotient, or a product) the group of automorphisms need no longer be
transitive, and $B_{k}$ need not be a polynomial: in general only
$B_{k}=(\text{polynomial in }k)+O(k^{-\infty})$ follows. For a polarized flat torus, however,
the expansion is exact: the translations act transitively and lift to $L^{k}$, so $B_{k}$ is
constant, and $h^{0}(M,L^{k})=k^{n}\chi(L)$ by Riemann--Roch and Kodaira vanishing, so that
$B_{k}=h^{0}(M,L^{k})/\Vol_{g}(M)=(k/\pi)^{n}$. Indeed each $a_{j}$ is a
universal polynomial in the curvature and its covariant derivatives
\cite{Lu00}, so the $a_{j}$ of such an $(M,g;L,h)$ coincide with those of its
model; hence $a_{j}=0$ for $j>n$ there as well, and $X$ is Bergman
logarithmically flat. This is \cite[Cor.~1.13]{EXX}.

The argument above gives an elementary proof of the statement on
logarithmic flatness in \cite[Cor.~1.13]{EXX} for these models: for
$\CP^{r}$ it uses only Riemann--Roch, and for $\BB^{l}(q)$ and $\mathbb C^{n}$ the explicit kernels. Local rigidity of Bergman
logarithmic flatness near the invariant metrics of compact Hermitian
symmetric spaces is proved in \cite{LS26}. All their coefficients $a_{j}$ with $j>n$ vanish, so the
vanishing pattern does not distinguish the models.
\end{remark}

\begin{corollary}[Local rigidity]\label{cf:cor:rigid}
Let $M=\CP^{n}$, $L=\mathcal O(1)$, and let $\omega\in[\omega_{\FS}]$ be real
analytic. Suppose that the logarithmic term of $\Pi$ vanishes to infinite
order, in the sense of Definition~\ref{cf:def:log}, over a non-empty open subset of $\CP^{n}$.
Then $a^{L}_{j}=0$ on $\CP^{n}$ for every $j\ge n+1$. Consequently, if $n=1$, or if $\omega$ is
the curvature of a metric $h$ on $\mathcal O(1)$ with
$\|h/h_{\FS}-1\|_{\mathscr{C}^{2n+4}}<\varepsilon_{n}$, where $\varepsilon_{n}$ is the constant
of \cite[Thm.~1.2]{LT}, then
$\omega=F^{*}\omega_{\FS}$ for an automorphism $F$ of $\CP^{n}$.
\end{corollary}
\begin{proof}
By the implication (i)$\Rightarrow$(iii) of Theorem~\ref{thm:logS}, $a^{L}_{j}=0$ for
$j\ge n+1$ on the given open set. The $a^{L}_{j}$ are universal polynomials in the curvature of
$g$ and its covariant derivatives \cite{Lu00}, hence real analytic, and $\CP^{n}$ is connected;
so $a^{L}_{j}=0$ on $\CP^{n}$ for every $j\ge n+1$. The theorems \cite[Thms.~1.1, 1.2]{LT}
assume the vanishing of $v$ in \cite[(2.4)]{LT}, and their proofs use it only through
\cite[Thm.~2.3]{LT}, that is, through $a^{L}_{n+1}=0$: the proof of \cite[Thm.~1.1]{LT} in
\cite[\S6]{LT} starts from $a_{2}=0$, and that of \cite[Thm.~1.2]{LT} in \cite[\S7]{LT} from
$a_{n+1}=0$ \cite[(7.3)]{LT}. These proofs therefore apply.
\end{proof}
Corollary~\ref{cf:cor:rigid} is the theorem of Lu and Tian with their hypothesis replaced by the
vanishing of the logarithmic term to infinite order over a non-empty open set, for real analytic
$\omega$. Its proof uses only the implication (i)$\Rightarrow$(iii) of Theorem~\ref{thm:logS},
not the converse.

\section{Compact Sasakian manifolds and dimension three}\label{sec:compact}
\subsection{Compact Sasakian manifolds}\label{ssec:sa:compact}
\begin{corollary}[Vanishing integrals]\label{cor:sa:integrals}
$\int_{X}a_{j}(g^{T})\,dV_{\eta}=0$ for every $j\ge n+1$.
\end{corollary}
\begin{proof}
Let $X$ be quasi-regular; multiplying $\eta$ by a constant, which multiplies
$a_{j}(g^{T})\,dV_{\eta}$ by a constant, we may assume that $T$ generates a circle
action of period $2\pi$. Integrating \eqref{eq:sa:diag} gives
$\operatorname{Tr}\big(\chi(-\sqrt{-1}T/k)\Pi\big)\sim\sum_{j}m_{n-j}(\chi)\big(\int_{X}s_{j}\,dV_{\eta}\big)k^{n+1-j}$.
The left side is $\sum_{m}\chi(m/k)\,d_{m}$, $d_{m}=\dim H^{0}_{b,m}(X)$. For $m$
large, $d_{m}=\sum_{i\le n}c_{i}(m)m^{i}$ with $c_{i}$ periodic of some period $N$,
by Kawasaki's Riemann--Roch theorem \cite{Kaw79} and the vanishing of higher cohomology \cite{Bai57} on
the quotient orbifold; see \cite{CHT19,CHT25} for the formulation on $X$.
Since $\operatorname{supp}\chi\subset(0,\infty)$ is compact, only the $m\ge ck$ contribute,
for some $c>0$, so we may replace $d_{m}$ by this quasi-polynomial and expand
each $c_{i}$ in the characters $m\mapsto e^{2\pi\sqrt{-1}rm/N}$, $0\le r<N$. By the
Poisson summation formula \cite[Thm.~7.2.1]{Hoe03},
\[
 \sum_{m\in\mathbb Z}\chi(m/k)\,m^{i}e^{2\pi\sqrt{-1}rm/N}
 =k^{i+1}\sum_{\ell\in\mathbb Z}\int\chi(t)\,t^{i}e^{2\pi\sqrt{-1}kt(r/N-\ell)}\,dt,
\]
and integration by parts in the terms with $\ell\ne r/N$ shows that this is
$m_{i}(\chi)k^{i+1}+O(k^{-\infty})$ if $r=0$, and $O(k^{-\infty})$ otherwise. So the
trace has no terms $k^{\ell}$ with $\ell\le0$, and the claim follows from
$m_{n-j}(\chi)>0$ and Lemma~\ref{lem:sa:transverse}. In general, with $\mathfrak t^{+}$ as
in Lemma~\ref{lem:sa:cone}, the integral for the structure of $T'\in\mathfrak t^{+}$
depends continuously on $T'$ and vanishes on a dense set; indeed
$a_{j}(g^{T})$ is a polynomial in the curvature of the Tanaka--Webster
connection of $\eta_{T'}$ on $\ker\eta$ and its covariant derivatives, which depend
smoothly on $T'$.
\end{proof}
\begin{remark}[The index character]\label{rem:sa:index}
The vanishing can also be read off from the index character of Martelli, Sparks and Yau
\cite{MSY08}, without approximating $T$. The cone over $X$ is an affine variety with an isolated
singularity at the vertex \cite[\S2.1]{CS18}. An eigenfunction of $-\sqrt{-1}T$ on the CR functions
with eigenvalue $\lambda$ extends to the cone as a holomorphic function, homogeneous of degree
$\lambda$, so the multiplicity of $\lambda$ is $d_{\lambda}=\sum_{\alpha(T)=\lambda}\dim H_{\alpha}$, where
$H=\bigoplus_{\alpha}H_{\alpha}$ is the weight decomposition of the coordinate ring under the torus
$\mathbb T$ of Lemma~\ref{lem:sa:cone}; thus $\sum_{\lambda}d_{\lambda}e^{-t\lambda}$ is the index
character $F(T,t)$ of \cite[Def.~4.1]{CS18}. By the proof of \cite[Thm.~4.10]{CS18}, $F(T,t)$ is a
finite sum of terms $\pm e^{-tb}\prod_{i=1}^{N}(1-e^{-tc_{i}})^{-1}$ with $b\in\mathbb R$ and $c_{i}>0$.
Hence $\operatorname{Tr}\big(\chi(-\sqrt{-1}T/k)\Pi\big)=\sum_{\lambda}\chi(\lambda/k)\,d_{\lambda}$ is a finite
sum of lattice sums $\pm\sum_{m\in\mathbb Z_{\ge0}^{N}}\chi\big((b+c\cdot m)/k\big)$. Since $\chi$
vanishes near $0$, the Euler--Maclaurin formula shows that each of them is, modulo $O(k^{-\infty})$, a
polynomial in $k$ without constant term, and comparing with the integral of \eqref{eq:sa:diag} gives
the claim. The proof above approximates $T$ by quasi-regular Reeb fields instead, as \cite{CS18} do
for the Donaldson--Futaki invariant.
\end{remark}

\begin{corollary}[Compact Sasakian manifolds]\label{cor:sa:compact}
The logarithmic term of $\Pi$ vanishes if and only if every
$a_{n+1+p}(g^{T})$ is constant, if and only if no $a_{n+1+p}(g^{T})$ changes
sign. If the Sasakian structure is real analytic, the vanishing over a
non-empty open set implies the vanishing over $X$.
\end{corollary}
\begin{proof}
A function that does not change sign and has integral zero vanishes, by
Corollary~\ref{cor:sa:integrals}; apply Theorem~\ref{thm:logS}. In the real
analytic case the functions $a_{j}(g^{T})$ are real analytic.
\end{proof}
This is the Sasakian form of \cite[Cor.~1.12]{EXX}.

\subsection{Dimension three}\label{ssec:sa:three}
\begin{lemma}[The first coefficients for $n=1$]\label{lem:sa:n1}
For $n=1$, with $\Delta=g^{1\bar1}\partial_{z}\partial_{\bar z}$,
$|\nabla f|^{2}=g^{1\bar1}|\partial_{z}f|^{2}$ and
$r_{g}=-g^{1\bar1}\partial_{z}\partial_{\bar z}\log g_{1\bar1}$, the scalar curvature in the
normalization of \cite{LT}, so that $r_{g}=2$ on $\CP^{1}$ for $\omega_{g}=\pi c_{1}(\mathcal O(1))$,
\begin{equation}\label{eq:sa:n1}
 a_{1}=\tfrac12r_{g},\qquad a_{2}=\tfrac13\Delta r_{g},\qquad
 a_{3}=\tfrac18\Delta^{2}r_{g}-\tfrac{5}{24}r_{g}\Delta r_{g}-\tfrac{5}{24}|\nabla r_{g}|^{2}
 =\tfrac18\Delta^{2}r_{g}-\tfrac{5}{48}\Delta(r_{g}^{2}).
\end{equation}
\end{lemma}
\begin{proof}
By Lu's theorem \cite[Thm.~1.1]{Lu00}, with $\rho$ the scalar curvature
normalized by $a_{1}=\frac12\rho$, so that $\rho=r_{g}$,
\begin{align*}
 a_{2}&=\tfrac13\Delta\rho+\tfrac1{24}\big(|R|^{2}-4|\Ric|^{2}+3\rho^{2}\big),\\
 a_{3}&=\tfrac18\Delta\Delta\rho+\tfrac1{24}\operatorname{div}\operatorname{div}(R,\Ric)
 -\tfrac16\operatorname{div}\operatorname{div}(\rho\Ric)
 +\tfrac1{48}\Delta\big(|R|^{2}-4|\Ric|^{2}+8\rho^{2}\big)\\
 &\quad+\tfrac1{48}\rho\big(\rho^{2}-4|\Ric|^{2}+|R|^{2}\big)
 +\tfrac1{24}\big(\sigma_{3}(\Ric)-\Ric(R,R)-R(\Ric,\Ric)\big).
\end{align*}
For $n=1$ the curvature tensor has a single component, $\Ric=\rho\,g$ and
$|R|^{2}=|\Ric|^{2}=\rho^{2}$. Hence $a_{2}=\frac13\Delta\rho$. By the definitions
\cite[(5.5)]{Lu00}, $\operatorname{div}\operatorname{div}(R,\Ric)$ is a double divergence of the
tensor $(R,\Ric)_{k\bar l}=\sum_{i,j}R_{i\bar jk\bar l}R_{\bar ji}$, which is linear in $R$, and
$R(\Ric,\Ric)$ is the contraction of the same tensor with $\Ric$, while $\Ric(R,R)$ is quadratic
in $R$. For $n=1$, $(R,\Ric)=\epsilon\rho^{2}g$ with a single sign $\epsilon=\pm1$ that depends on
the sign convention for the curvature tensor; then $R(\Ric,\Ric)=\epsilon\rho^{3}$ and
$\operatorname{div}\operatorname{div}(R,\Ric)=\epsilon\,\Delta(\rho^{2})$, while
$\sigma_{3}(\Ric)=\rho^{3}$ and $\Ric(R,R)=\rho^{3}$. Moreover
$\operatorname{div}\operatorname{div}(\rho\Ric)=\Delta(\rho^{2})$. Therefore
\[
 a_{3}=\tfrac18\Delta^{2}\rho+c\,\Delta(\rho^{2})+e\,\rho^{3},\qquad
 c=\tfrac{\epsilon}{24}-\tfrac16+\tfrac5{48},\qquad
 e=-\tfrac{1}{24}+\tfrac{1}{24}(1-1-\epsilon)=-\tfrac{1+\epsilon}{24}.
\]
The sign is determined by a model: on $\CP^{1}$ the curvature is a non-zero
constant and $a_{3}=0$ (Remark~\ref{rem:models}), so $e=0$, that is,
$\epsilon=-1$. Hence $c=-\frac5{48}$, and \eqref{eq:sa:n1} follows from
$\Delta(\rho^{2})=2\rho\Delta\rho+2|\nabla\rho|^{2}$.
\end{proof}
\begin{remark}\label{rem:sa:n1check}
The values of $a_{1}$ and $a_{2}$ agree with \cite[Thm.~2.2]{LT}; in the
normalization of \cite{MM07}, where $\omega=c_{1}(L)$, one has $r_{g}=r^{X}/(4\pi)$, and
$a_{1}=\frac12r_{g}$ is $b_{1}=r^{X}/(8\pi)$ of \cite[(4.1.8)]{MM07}. The coefficient $\frac18$
of $\Delta^{2}r_{g}$ is $k/(k+1)!$ for $k=3$ in \cite[Thm.~3.1]{LT}. As an independent
check of the sign $\epsilon$, for the radial weights $e^{-k(|z|^{2}+\varepsilon|z|^{4})}$ on $\mathbb C$ the
monomials are orthogonal and $B_{k}$ is an explicit series; its numerical
evaluation at points with $\nabla r_{g}\ne0$ gives the coefficients of
$r_{g}\Delta r_{g}$ and of $|\nabla r_{g}|^{2}$ as $-\frac5{24}$ to six digits, and excludes the
value $-\frac1{24}$ that corresponds to $\epsilon=+1$.
\end{remark}

\begin{corollary}[Dimension three]\label{cor:sa:three}
Let $n=1$ and let $W\subset X$ be a connected open set. The logarithmic term
of $\Pi$ vanishes to infinite order over $W$ if and only if
$a_{2}(g^{T})=a_{3}(g^{T})=0$ on $W$, if and only if the transverse curvature
is constant on $W$; then $W$ is locally isomorphic to a Sasakian space form,
and in particular locally CR spherical. Consequently, for compact $X$ the
logarithmic term vanishes if and only if the transverse curvature is
constant, and it never vanishes if $X$ is irregular.
\end{corollary}
\begin{proof}
By \eqref{eq:sa:n1}, $a_{2}=0$ on $W$ gives $\Delta r_{g^{T}}=0$, hence
$a_{3}=-\frac{5}{24}|\nabla r_{g^{T}}|^{2}$ there, and $a_{3}=0$ makes $r_{g^{T}}$ constant. Only
the fact that the coefficient of $|\nabla r_{g^{T}}|^{2}$ in $a_{3}$ is not zero is used
here; it does not depend on the normalization, since a change of scale
multiplies the three terms of $a_{3}$ by the same positive factor, and a change
of sign of $\Delta$ affects only the term $r\Delta r$, which vanishes on $W$. For
compact $X$ the condition $a_{2}=0$ alone suffices, since a harmonic basic
function on a compact Sasakian manifold is constant; the coefficient $a_{3}$ is
needed for open sets $W$, on which harmonic functions need not be
constant, and gives the local statement, which is finer than what the
first logarithmic coefficient yields \cite{Hir93}. A
Sasakian structure of constant transverse curvature is locally a Sasakian
space form \cite{Tan69}, whose transverse geometry is that of $\CP^{1}$,
$\mathbb C$ or the disc; there $a_{j}=0$ for $j\ge2$ by Remark~\ref{rem:models}
and locality, for every scale since $a_{j}(cg)=c^{-j}a_{j}(g)$ for constants
$c>0$, and Theorem~\ref{thm:logS} applies. In dimension three such
structures have vanishing torsion and constant Webster curvature, hence
vanishing Cartan tensor. For compact $X$ take $W=X$. Finally, a compact Sasakian $3$-manifold of
constant transverse curvature is $\Gamma\backslash\widetilde X$ with $\widetilde X$ the
standard structure on $S^{3}$, the Heisenberg group or
$\widetilde{SL}(2,\mathbb R)$ \cite{Tan69}, and the Reeb flow is induced by the
vertical flow of the model: the Hopf action, the action of the center,
which meets every lattice in a lattice, or the fibers of the Seifert
fibration \cite{Sco83}. Its orbits are closed, so $X$ is quasi-regular; see
also \cite{Gei97}, \cite[\S3]{Bel01}.
\end{proof}

\begin{remark}[Comparison]\label{rem:sa:compare}
The compact part of Corollary~\ref{cor:sa:three}, and its consequence for
irregular $X$, also follow from Hirachi's formula
$\psi|_{M}=\frac{1}{24\pi^{2}}(\Delta_{b}R-2\operatorname{Im}A^{11}{}_{,11})$ for the first logarithmic
coefficient \cite{Hir93}, since the torsion vanishes for $\theta=\eta$. The local
part uses the second coefficient. It differs from Hirachi's
Corollary~3 in \cite{Hir93} (any contact form, compact boundary with a
transversal symmetry, vanishing to third order, conclusion: spherical) and
from Ebenfelt's \cite[Cor.~1.3]{Ebe18} (the invariant volume element,
compact boundary with a transversal symmetry, vanishing to second order,
conclusion: spherical): here the volume element is $dV_{\eta}$, the
hypothesis is local, and the conclusion, constant $R$, is stronger than
sphericity, which for $\theta=\eta$ means $R_{;\bar1\bar1}=0$ \cite[Thm.~4.2]{Ebe18}.
Sphericity does not suffice. On the standard CR sphere $S^{3}\subset\mathbb C^{2}$,
take the Reeb field $a_{1}\partial_{\theta_{1}}+a_{2}\partial_{\theta_{2}}$ of the Sasaki
cone with $a_{1}\ne a_{2}$. If $a_{1}/a_{2}\notin\mathbb Q$, the structure is irregular
and $R$ is not constant by Corollary~\ref{cor:sa:three}. If $a_{1}/a_{2}=p/q$ with
$p\ne q$ coprime, the leaf space is the spindle $S^{2}(p,q)$, a bad orbifold,
which carries no metric of constant curvature \cite[\S2]{Sco83}. In both
cases the first logarithmic coefficient is not identically zero, although
the CR structure is spherical. Hirachi notes that
a localized statement requires a choice of volume element \cite{Hir93}; the
rigidity statement of Hanges concerns the second, non-logarithmic
coefficient, for the volume element $dx\,dy\,dt$ in rigid coordinates \cite[(2.2), Rem.~4.4]{Han90}.
\end{remark}

\section{$(n,0)$-forms on complete Sasakian manifolds}\label{tf:sec:setting}
This section proves Theorem~\ref{thm:introC}. Compared with the case of functions, the only change of
method is the semi-classical cut-off (Section~\ref{tf:ssec:cutoff}). The local structure of the
Szeg\H o projection is \cite[Thm.~1.3]{HMW22}, and the remaining steps are those of
Sections~\ref{sec:nf}--\ref{sec:log} and Appendix~\ref{app:nf}.

\subsection{Setting}\label{tf:ssec:setting}
Let $(X,\eta,T)$ be a Sasakian manifold of dimension $2n+1$ whose Sasakian metric is complete; $X$ need
not be compact. This is the setting of \cite[Assumption~1.1, Thm.~1.3]{HMW22}: the flow $\varphi_{s}$ of $T$
is an $\mathbb R$-action which preserves $\eta$ and the CR structure, and the Levi metric of
\cite[(1.4)]{HMW22} is the Sasakian metric. We normalize as in Section~\ref{sec:sasaki}: in BRT
coordinates $(z,\vartheta)$, $T=\partial_{\vartheta}$ and $d\eta=2\omega_{g}$ with
$\omega_{g}=\frac{\sqrt{-1}}{2}\partial\overline\partial\varphi$.

Let $L^{2}_{(n,0)}(X)$ be the space of $(n,0)$-forms $\alpha$ with $\langle\alpha,\alpha\rangle<\infty$, where
\begin{equation}\label{tf:eq:pairing}
 \langle\alpha,\beta\rangle=\int_{X}2^{-n}(\sqrt{-1})^{n^{2}}\,\alpha\wedge\overline{\beta}\wedge\eta .
\end{equation}
This pairing involves no metric; for every $\mathbb R$-invariant Hermitian metric as in
\cite[(1.5)]{HMW22}, the $L^{2}$ structure of \cite{HMW22} on $(n,0)$-forms is a constant multiple of it
(Lemma~\ref{tf:lem:vol}). The \emph{Hardy space of $(n,0)$-forms} is
\[
 H^{n,0}_{b}(X)=\{\alpha\in L^{2}_{(n,0)}(X):\overline\partial_{b}\alpha=0\},
\]
with $\overline\partial_{b}\alpha$ in the sense of distributions; it is a closed subspace, and
$S:L^{2}_{(n,0)}(X)\to H^{n,0}_{b}(X)$ denotes the orthogonal projection, the Szeg\H o projection on
$(n,0)$-forms.

On a BRT chart we write $(n,0)$-forms as $u\,e$ with $e=dz_{1}\wedge\dots\wedge dz_{n}$; $d\lambda(z)$ is the
Lebesgue measure in $z$, and $\det\varphi''=\det(\partial_{z_{j}}\partial_{\overline z_{k}}\varphi)$.

\begin{lemma}[Volume forms]\label{tf:lem:vol}
On a BRT chart, $dV_{\eta}=\det\varphi''\,d\lambda(z)\wedge d\vartheta$, and
$\langle u\,e,v\,e\rangle=\int u\overline v\,dV_{e}$ with $dV_{e}:=d\lambda(z)\wedge d\vartheta=|e|^{2}dV_{\eta}$, where
$|e|^{2}=(\det\varphi'')^{-1}$ is the pointwise norm defined by \eqref{tf:eq:pairing}; in the regular case it
is the norm on $K_{M}$ induced by $g$. For an $\mathbb R$-invariant Hermitian metric $g_{X}$ as in
\cite[(1.5)]{HMW22}, with volume form $dv_{X}$, $|e|^{2}_{g_{X}}dv_{X}$ is a constant multiple of $dV_{e}$.
\end{lemma}

\begin{proof}
Write $\eta=d\vartheta+\beta$ with $\beta$ a $1$-form in $z$. Then $\beta\wedge(d\eta)^{n}=0$ by degree and
$(d\eta)^{n}/n!=2^{n}\det\varphi''\,d\lambda$, so $dV_{\eta}=\det\varphi''\,d\lambda\wedge d\vartheta$; and
$(\sqrt{-1})^{n^{2}}e\wedge\overline e=2^{n}d\lambda$, so $\langle u\,e,v\,e\rangle=\int u\overline v\,d\lambda\,d\vartheta$. The
last assertion is \cite[(4.56)--(4.59)]{HMW22}.
\end{proof}

\subsection{The semi-classical cut-off}\label{tf:ssec:cutoff}
Let $U(s)\alpha=\varphi_{s}^{*}\alpha$ and $A=-\sqrt{-1}\mathcal L_{T}$.

\begin{lemma}\label{tf:lem:sa}
$U(s)$ is a strongly continuous one-parameter group of unitary operators on $L^{2}_{(n,0)}(X)$, so
$U(s)=e^{\sqrt{-1}sA}$ with $A$ self-adjoint. Moreover $U(s)S=SU(s)$; hence $f(A)S=Sf(A)$ for every
bounded Borel function $f$, and the kernel of $f(A)S$ is invariant under $\varphi_{s}\times\varphi_{s}$.
\end{lemma}

\begin{proof}
The flow preserves $\eta$ and the CR structure, hence the pairing \eqref{tf:eq:pairing}, so $U(s)$ is
unitary; it is strongly continuous on the dense subspace of compactly supported smooth forms, hence
everywhere, and Stone's theorem gives $A$. Since the flow is CR, $U(s)$ maps $H^{n,0}_{b}(X)$ into
itself for every $s$, hence onto itself. A unitary operator which maps a closed subspace onto itself
maps its orthogonal complement onto itself, so it commutes with $S$. The remaining assertions follow
from the spectral theorem.
\end{proof}

For $\chi\in\mathscr{C}^{\infty}_{c}((0,\infty))$ we use the cut-off $\chi(k^{-1}A)S=S\chi(k^{-1}A)S$. On a
compact manifold the Toeplitz operator $SAS=S(-\sqrt{-1}\mathcal L_{T})S$ is self-adjoint and, since $A$
commutes with $S$, $\chi(k^{-1}SAS)=\chi(k^{-1}A)S$, as for functions in \eqref{eq:sa:chiT}.
On a non-compact manifold we do not use the Toeplitz operator, whose self-adjointness would have to be
established: the operator $A$ is self-adjoint on all of $L^{2}_{(n,0)}(X)$, and this is the only global
property of $X$, besides \cite[Thm.~1.3]{HMW22}, which enters below. It uses that the Reeb flow
preserves the $L^{2}$ structure; for a transversal CR vector field whose flow does not, neither the
unitarity nor the commutation with $S$ is available.

\subsection{The Szeg\H o kernel in a frame}\label{tf:ssec:frame}
On a BRT chart $D\Subset X$ let $K$ be the kernel of $S$ in the frame $e$ with respect to $dV_{e}$, that is,
$S(u\,e)(x)=\big(\int K(x,y)u(y)\,dV_{e}(y)\big)e(x)$ for $u\in\mathscr C^{\infty}_{c}(D)$.

\begin{proposition}[Local structure]\label{tf:prop:H}
\begin{enumerate}
\item The kernel of $S$ is smooth on $(X\times X)\setminus\operatorname{diag}$.
\item On every BRT chart $D\Subset X$, $K$ satisfies \ref{K:rep}--\ref{K:sym} with $dV=dV_{e}$, and the
      leading coefficient in \ref{K:rep} is $\sigma_{0}(x,x)=\det\varphi''(z)/(2\pi^{n+1})$.
\item Consequently $K\equiv I(\sigma)$ with $\sigma\sim\sum_{j}\sigma_{j}(z,w)t^{n-j}$ independent of the Reeb
      variables, unique modulo terms flat at $z=w$; $\partial_{\overline z}\sigma_{j}$ and $\partial_{w}\sigma_{j}$
      are flat, and $\sigma_{j}$ is determined modulo flat terms by $\sigma_{j}(z,z)$.
\item The symbol $\sigma$ is determined, modulo terms flat at $z=w$, by the CR structure and $\eta$ on $D$.
\end{enumerate}
\end{proposition}

\begin{proof}
(1) and (2). Two points of the connected manifold $X$ lie in a relatively compact coordinate patch, so
it suffices to work on such a patch $D$, which for (2) we take to be a BRT chart. By
\cite[Thm.~1.3]{HMW22}, in the trivialization $e$ the kernel of $S$ is
$\int_{0}^{\infty}e^{\sqrt{-1}t\varphi_{-}}s_{-}\,dt$ modulo $\mathscr C^{\infty}(D\times D)$, where $s_{-}$ is a
classical symbol of order $n$ and $\varphi_{-}$ satisfies \cite[(1.9)]{HMW22}; in particular
$\varphi_{-}(x,y)\ne0$ for $x\ne y$, and integration by parts in $t$ gives (1). The phase is that of
\cite{HM17JDG} constructed on $D$ \cite[p.~5]{HMW22}, so it satisfies \textup{($\alpha$)--($\gamma$)} of
Theorem~\ref{thm:BdMS} \cite[Thms.~3.4, 4.4]{HM17JDG}, and Proposition~\ref{nf:prop:phase} gives
$K\equiv I(\sigma')$ with the BRT phase. By \cite[Thm.~4.5, (4.2)]{HM17JDG}, the leading coefficient on the
diagonal is a normalizing constant times the ratio of $dV_{\eta}$ to the volume form with respect to
which the kernel is taken; the constant is $(2\pi^{n+1})^{-1}$, as for $dV_{\eta}$ \cite[(1.1.4)]{CHH25},
and the ratio is $dV_{\eta}/dV_{e}=\det\varphi''$ by Lemma~\ref{tf:lem:vol}. This is \ref{K:rep}. Since
$\mathcal L_{T}e=0$ and $dV_{e}$ is invariant, $U(s)(u\,e)=(u\circ\varphi_{s})\,e$, and Lemma~\ref{tf:lem:sa}
gives \ref{K:inv}. Since $\overline\partial_{b}(u\,e)=(\overline\partial_{b}u)\wedge e$ \cite[(2.18)]{HMW22} and $S$
maps into $H^{n,0}_{b}(X)$, $\overline Z_{l,x}K=0$, which is \ref{K:CR}. Finally $S$ is self-adjoint for
\eqref{tf:eq:pairing}, which in the frame $e$ is $\int u\overline v\,dV_{e}$ (Lemma~\ref{tf:lem:vol}); this is
\ref{K:sym}.

(3) Proposition~\ref{nf:prop:symbol}.

(4) In the proof of \cite[Thm.~1.3]{HMW22}, which is that of \cite[Thm.~4.35]{HMW22} with
\cite[Thm.~4.15]{HMW22} in place of \cite[Thm.~4.13]{HMW22}, the phase and the symbol are those of the
operator $\widetilde S$ of \cite[Thm.~4.30, (4.89)]{HMW22}, which is constructed on $D$ from the CR
structure and the volume form there, see \cite[Thms.~4.31, 4.34]{HMW22}. Two charts with the same data
therefore give phases and symbols which agree modulo flat terms, and the uniqueness in (3) gives the
claim.
\end{proof}

\begin{remark}[Change of frame]\label{tf:rem:frame}
For another holomorphic frame $e'=f(z)e$ of the canonical bundle of the leaf space, the scalar kernel is
multiplied by $f(z)^{-1}\overline{f(w)}^{-1}$, up to the factor relating $dV_{e}$ and $dV_{e'}$, a non-vanishing
almost sesquiholomorphic factor independent of the Reeb variables. So the class of Reeb-independent
symbols, the vanishing of the logarithmic term, and the products $|e|^{2}\sigma_{j}(x,x)$ do not depend on
the frame.
\end{remark}

\subsection{Localization of the cut-off}\label{tf:ssec:loc}
Write $\chi(\tau)=\frac{1}{2\pi}\int\widehat\chi(s)e^{\sqrt{-1}s\tau}ds$, fix
$\zeta\in\mathscr C^{\infty}_{c}((-\varepsilon,\varepsilon))$ with $\zeta=1$ near $0$, where $\varepsilon>0$ is small,
and split
\begin{align*}
 \chi(k^{-1}A)S&=L_{k}+R_{k}S, & L_{k}&=\frac{1}{2\pi}\int k\widehat\chi(ks)\,\zeta(s)\,U(s)S\,ds,\\
 R_{k}&=r_{k}(A), & r_{k}(\tau)&=\frac{1}{2\pi}\int k\widehat\chi(ks)\,(1-\zeta(s))\,e^{\sqrt{-1}s\tau}ds .
\end{align*}

\begin{lemma}\label{tf:lem:rk}
For all $N$ and $L$, $\sup_{\tau\in\mathbb R}(1+|\tau|)^{L}|r_{k}(\tau)|=O(k^{-N})$.
\end{lemma}

\begin{proof}
Integrate by parts in $s$: $\widehat\chi$ is a Schwartz function, and $|s|$ is bounded below on the
support of $1-\zeta$.
\end{proof}

\begin{proposition}[Local part]\label{tf:prop:local}
Let $D\Subset X$ be a BRT chart and $\sigma$ as in Proposition~\ref{tf:prop:H}. On compact subsets of
$D\times D$, in the frame $e$,
\[
 L_{k}(x,y)=\int_{0}^{\infty}e^{\sqrt{-1}t\phi(x,y)}\chi(t/k)\,\sigma(z,w,t)\,dt+O(k^{-\infty})
\]
in the $\mathscr C^{\infty}$ topology, and on compact subsets of $(X\times X)\setminus\operatorname{diag}$ the kernel
of $L_{k}$ is $O(k^{-\infty})$.
\end{proposition}

\begin{proof}
For $|s|<\varepsilon$ and $x$ in a compact subset of $D$, $\varphi_{s}x\in D$ and the kernel of $U(s)S$ is
$K(\varphi_{s}x,y)$. Write $K=I(\sigma)+F$ with $F$ smooth. By Taylor's formula in $s$, the contribution of
$F$ is $\sum_{j<N}k^{-j}(-\sqrt{-1})^{j}\chi^{(j)}(0)\,\partial_{\vartheta}^{j}F/j!+O(k^{-N})$, and
$\chi^{(j)}(0)=0$. For the oscillatory part, $\phi(\varphi_{s}x,y)=\phi(x,y)+s$ by (P1), and $\sigma$ does not
depend on the Reeb variables, so the integral in $s$ gives
$\int_{0}^{\infty}e^{\sqrt{-1}t\phi}(\chi(t/k)-r_{k}(t))\sigma\,dt$; since $\operatorname{Im}\phi\ge0$,
Lemma~\ref{tf:lem:rk} shows that the term with $r_{k}$ is $O(k^{-\infty})$. Off the diagonal the same
Taylor argument applies to the kernel of $S$, which is smooth there by Proposition~\ref{tf:prop:H}(1).
\end{proof}

For $W\Subset X$ open and $\alpha,\beta\in\mathscr C^{\infty}_{c}(W,T^{*n,0}X)$, let $\mu_{\alpha}$ and
$\mu_{\alpha,\beta}$ be the spectral measures of $S\alpha$ and of the pair $(S\alpha,S\beta)$ for $A$, so that
$\langle f(A)S\alpha,S\beta\rangle=\int f\,d\mu_{\alpha,\beta}$ and $\mu_{\alpha}=\mu_{\alpha,\alpha}$. Fix Sobolev norms
$\|\cdot\|_{H^{m}}$ on forms supported in $W$.

\begin{lemma}[Spectral measures]\label{tf:lem:weyl}
For every $m\ge0$ there is $C_{m}$ such that
$\mu_{\alpha}([\lambda,\lambda+1])\le C_{m}(1+|\lambda|)^{n+2m}\|\alpha\|^{2}_{H^{-m}}$ for all $\lambda\in\mathbb R$ and all
$\alpha\in\mathscr C^{\infty}_{c}(W,T^{*n,0}X)$.
\end{lemma}

\begin{proof}
Choose $\beta_{0}\in\mathscr C^{\infty}_{c}((-\varepsilon,\varepsilon))$ with $\beta_{0}\ge0$, $\int\beta_{0}=1$ and
$\varepsilon\le\frac12$; then $|\widehat\beta_{0}|\ge\frac12$ on $[-1,1]$, and $\varrho=|\widehat\beta_{0}|^{2}$ is the
Fourier transform of a function supported in $(-2\varepsilon,2\varepsilon)$. With $\lambda'=\lambda+\frac12$ we have
$\varrho(A-\lambda')\ge\frac14\,\mathbf 1_{[\lambda,\lambda+1]}(A)$, hence
$\mu_{\alpha}([\lambda,\lambda+1])\le4\langle\varrho(A-\lambda')S\alpha,\alpha\rangle$, and $\varrho(A-\lambda')S$ is an integral
of $U(s)S$ over $|s|<2\varepsilon$. As in Proposition~\ref{tf:prop:local}, covering $\overline W$ by finitely
many BRT charts, its kernel on $W\times W$ is $\int_{0}^{\infty}e^{\sqrt{-1}t\phi}\varrho(t-\lambda')\sigma\,dt$ near
the diagonal and is smooth away from it, modulo $O((1+|\lambda|)^{-\infty})$; so its derivatives of order
$\le2m$ are $O((1+|\lambda|)^{n+2m})$, uniformly on $W\times W$. A kernel of class $\mathscr C^{2m}$ on
$\overline W\times\overline W$ defines a form bounded by its $\mathscr C^{2m}$ norm times
$\|\alpha\|_{H^{-m}}^{2}$.
\end{proof}

\begin{theorem}[Localization]\label{tf:thm:main}
Let $D\Subset X$ be a BRT chart and $\sigma$ as in Proposition~\ref{tf:prop:H}. For every
$\chi\in\mathscr{C}^{\infty}_{c}((0,\infty))$, on compact subsets of $D\times D$, in the frame $e$,
\[
 (\chi(k^{-1}A)S)(x,y)=\int_{0}^{\infty}e^{\sqrt{-1}t\phi(x,y)}\chi(t/k)\,\sigma(z,w,t)\,dt+O(k^{-\infty}),
\]
and on compact subsets of $(X\times X)\setminus\operatorname{diag}$ the kernel of $\chi(k^{-1}A)S$ is
$O(k^{-\infty})$.
\end{theorem}

\begin{proof}
By Proposition~\ref{tf:prop:local} it suffices to show that the kernel of $R_{k}S$ is $O(k^{-\infty})$ in
$\mathscr C^{\infty}(W\times W)$ for every open $W\Subset X$. Since $R_{k}$ commutes with $S$,
$\langle R_{k}S\alpha,\beta\rangle=\int r_{k}\,d\mu_{\alpha,\beta}$, and $|\mu_{\alpha,\beta}|(B)\le\mu_{\alpha}(B)^{1/2}\mu_{\beta}(B)^{1/2}$
for Borel sets $B$ by the Cauchy--Schwarz inequality. Splitting $\mathbb R$ into the intervals
$[\lambda,\lambda+1]$, $\lambda\in\mathbb Z$, Lemmas~\ref{tf:lem:rk} and~\ref{tf:lem:weyl} give
$|\langle R_{k}S\alpha,\beta\rangle|\le C_{N,m}k^{-N}\|\alpha\|_{H^{-m}}\|\beta\|_{H^{-m}}$ for all $N$ and $m$. Hence
$R_{k}S$ maps forms in $H^{-m}$ supported in $W$ to $H^{m}(W)$ with norm $O(k^{-\infty})$, for every $m$,
and its kernel is $O(k^{-\infty})$ in $\mathscr C^{\infty}(W\times W)$ by the Sobolev embedding theorem.
\end{proof}

\begin{remark}
The cut-off by $A$ is used only here: the part $R_{k}S$, which collects the returns of the Reeb flow,
is controlled through the spectral measures of the self-adjoint operator $A$, so that no compactness,
bounded geometry or control of the returns is needed for this step. For functions the same argument
applies wherever the Szeg\H o projection has the local structure of Proposition~\ref{tf:prop:H}.
\end{remark}

\subsection{Circle bundles}\label{tf:ssec:circle}
Let $(M,g)$ be a complete K\"ahler manifold of dimension $n$ and $(L,h)\to M$ a holomorphic line bundle
with $\omega_{g}=\frac{\sqrt{-1}}{2}\partial\overline\partial\varphi$ for $|e_{L}|^{2}_{h}=e^{-\varphi}$, and let
$\pi:X\to M$ be the circle bundle of $(L^{*},h^{*})$. It is Sasakian, its Sasakian metric is complete
\cite[Ex.~4.6]{HMW22}, and the Reeb flow is the circle action $r_{\theta}$. Let $E_{m}=L^{m}\otimes K_{M}$,
with the metric $h^{m}\otimes H$, where $H$ is the metric on $K_{M}$ induced by $g$, so that
$|dz_{1}\wedge\dots\wedge dz_{n}|^{2}_{H}=(\det\varphi'')^{-1}$. Let $H^{0}_{(2)}(M,E_{m})$ be the space of
holomorphic sections which are square integrable for $dV_{g}=\omega_{g}^{n}/n!$, and
$B^{E}_{m}=\sum_{i}|s_{i}|^{2}_{h^{m}\otimes H}$ for an orthonormal basis $(s_{i})$. Let
$H^{n,0}_{b,m}(X)=\{\alpha\in H^{n,0}_{b}(X):r_{\theta}^{*}\alpha=e^{\sqrt{-1}m\theta}\alpha\}$, and let $S_{m}$ be the
orthogonal projection onto it, with kernel $K_{m}$ in the frame $e$.

\begin{lemma}[Fourier components]\label{tf:lem:fourier}
In the coordinates $(z,\vartheta)$ of Lemma~\ref{lem:sa:BRT}, the map
$s=f\,e_{L}^{m}\otimes dz_{1}\wedge\dots\wedge dz_{n}\mapsto\alpha_{s}=f(z)e^{-m\varphi(z)/2}e^{\sqrt{-1}m\vartheta}\,e$ is an
isomorphism of $H^{0}_{(2)}(M,E_{m})$ onto $H^{n,0}_{b,m}(X)$ with $\|\alpha_{s}\|^{2}=2\pi\|s\|^{2}$. Consequently
\begin{equation}\label{tf:eq:KmB}
 K_{m}(x,x)=\frac{B^{E}_{m}(\pi(x))}{2\pi\,|e(x)|^{2}},\qquad |e|^{2}=(\det\varphi'')^{-1}.
\end{equation}
\end{lemma}

\begin{proof}
With $\overline Z_{l}=\partial_{\overline z_{l}}-\frac{\sqrt{-1}}{2}\varphi_{\overline z_{l}}\partial_{\vartheta}$,
\[
 \overline Z_{l}\big(f e^{-m\varphi/2}e^{\sqrt{-1}m\vartheta}\big)
 =(\partial_{\overline z_{l}}f)\,e^{-m\varphi/2}e^{\sqrt{-1}m\vartheta},
\]
so $\alpha_{s}$ is CR if and only if $f$ is holomorphic; conversely, a CR form $u\,e$ of weight $m$ has
$u=f_{1}(z)e^{\sqrt{-1}m\vartheta}$, and the same computation shows that $f_{1}e^{m\varphi/2}$ is holomorphic. The
definition does not depend on the frame of $L$ and the coordinates. By Lemma~\ref{tf:lem:vol},
$\|\alpha_{s}\|^{2}=\int|f|^{2}e^{-m\varphi}\,d\lambda\,d\vartheta=2\pi\int_{M}|f|^{2}e^{-m\varphi}\,d\lambda$, while
$|s|^{2}_{h^{m}\otimes H}\,dV_{g}=|f|^{2}e^{-m\varphi}\,d\lambda$. If $(s_{i})$ is an orthonormal basis of
$H^{0}_{(2)}(M,E_{m})$, then $(\alpha_{s_{i}}/\sqrt{2\pi})$ is one of $H^{n,0}_{b,m}(X)$, so
$K_{m}(x,x)=\frac{1}{2\pi}\sum_{i}|f_{i}|^{2}e^{-m\varphi}=\frac{1}{2\pi}\det\varphi''\,B^{E}_{m}$.
\end{proof}

\begin{proposition}[The expansion of $B^{E}_{m}$]\label{tf:prop:diag}
Let $\sigma$ be the symbol of Proposition~\ref{tf:prop:H}. Locally uniformly on $M$, with all derivatives,
\[
 B^{E}_{m}(z)=2\pi\,|e|^{2}\,\sigma(z,z,m)+O(m^{-\infty})
 \sim\pi^{-n}\sum_{j\ge0}a^{E}_{j}(z)\,m^{n-j},\qquad
 a^{E}_{j}:=2\pi^{n+1}|e|^{2}\sigma_{j}(z,z).
\]
If $M$ is compact, the $a^{E}_{j}$ are the coefficients of the Bergman kernel expansion of
$L^{m}\otimes K_{M}$ \cite[Thm.~4.1.1]{MM07}.
\end{proposition}

\begin{proof}
The form $e=\pi^{*}(dz_{1}\wedge\dots\wedge dz_{n})$ and $dV_{e}$ are invariant under the circle action, so
$K$ is defined on $\pi^{-1}(U)\times\pi^{-1}(U)$ for a coordinate ball $U\subset M$, and the kernel of
$U(\theta)S$ is $K(r_{\theta}x,y)$. By Lemma~\ref{tf:lem:sa},
$S_{m}=\frac{1}{2\pi}\int_{-\pi}^{\pi}e^{-\sqrt{-1}m\theta}U(\theta)S\,d\theta$, hence, as in \eqref{eq:fourier},
\[
 K_{m}(x,x)=\frac{1}{2\pi}\int_{-\pi}^{\pi}e^{-\sqrt{-1}m\theta}K(r_{\theta}x,x)\,d\theta ,
\]
and $B^{E}_{m}(\pi(x))=2\pi|e(x)|^{2}K_{m}(x,x)$ by \eqref{tf:eq:KmB}. Let
$\zeta\in\mathscr C^{\infty}_{c}((-\varepsilon,\varepsilon))$ equal $1$ near $0$. For $\theta\in[-\pi,\pi]$ outside a
neighborhood of $0$, $r_{\theta}x\ne x$ because the action is free, so $(1-\zeta(\theta))K(r_{\theta}x,x)$ is
smooth and $2\pi$-periodic in $\theta$, locally uniformly in $x$ with its derivatives
(Proposition~\ref{tf:prop:H}(1)), and its $m$-th Fourier coefficient is $O(m^{-\infty})$. For
$|\theta|<\varepsilon$, $\phi(r_{\theta}x,x)=\theta$ by (P2), so
$K(r_{\theta}x,x)=\int_{0}^{\infty}e^{\sqrt{-1}t\theta}\sigma(z,z,t)\,dt+F(r_{\theta}x,x)$ with $F$ smooth. The term
with $F$ contributes $O(m^{-\infty})$, and the oscillatory term contributes
\[
 \frac{1}{2\pi}\int_{0}^{\infty}\widehat\zeta(m-t)\,\sigma(z,z,t)\,dt=\sigma(z,z,m)+O(m^{-\infty}),
 \qquad \widehat\zeta(\tau)=\int e^{-\sqrt{-1}\tau\theta}\zeta(\theta)\,d\theta ,
\]
because $\frac{1}{2\pi}\int\widehat\zeta=\zeta(0)=1$, all moments of $\widehat\zeta$ of positive order vanish since
$\zeta=1$ near $0$, and Taylor's formula for $t\mapsto\sigma(z,z,t)$ at $t=m$, with the symbol estimates and
the decay of $\widehat\zeta$, gives the remainder. This proves the expansion, with
$a^{E}_{j}=2\pi^{n+1}|e|^{2}\sigma_{j}(z,z)$. If $M$ is compact, $B^{E}_{m}\sim(m/\pi)^{n}\sum_{j}a^{E}_{j}m^{-j}$
with the usual coefficients \cite[Thm.~4.1.1]{MM07}, and an asymptotic expansion in powers of $m$ is
unique.
\end{proof}

\begin{remark}[Complete bases]\label{tf:rem:bergman}
For a complete K\"ahler manifold whose K\"ahler form is the curvature of $L$, the expansion of the
Bergman function of $L^{m}\otimes K_{M}$ holds without further curvature or volume hypothesis. This is
due to Hsiao and Marinescu \cite[Thm.~1.7]{HM14}, who allow $L$ semi-positive, on the set where its
curvature is positive, and compute the first coefficients; see also \cite[Thm.~6.1.1]{MM07}.
Proposition~\ref{tf:prop:diag} recovers it from the Szeg\H o projection of the circle bundle. For $L^{m}$
itself, compare \cite[Thm.~1.6]{HMW22}, where curvature conditions are needed.
\end{remark}

\subsection{Identification of the coefficients and the logarithmic term}\label{tf:ssec:general}

\begin{theorem}[Identification]\label{tf:thm:id}
Let $D\Subset X$ be a BRT chart and $\sigma$ as in Proposition~\ref{tf:prop:H}. Then, on $D$ and for every $j$,
\[
 |e(x)|^{2}\sigma_{j}(x,x)=\frac{a^{E}_{j}(g^{T})(x)}{2\pi^{n+1}},
\]
where $a^{E}_{j}(\cdot)$ is the universal curvature polynomial giving the $j$-th coefficient of the Bergman
kernel expansion of $L^{k}\otimes K_{M}$ for a positive line bundle $L$ over a K\"ahler manifold $M$
\cite[Thm.~4.1.2]{MM07}, evaluated on $g^{T}$.
\end{theorem}

\begin{proof}
Fix $x_{0}\in D$. By Remark~\ref{tf:rem:frame} and the uniqueness in Proposition~\ref{tf:prop:H}(3), the
products $|e|^{2}\sigma_{j}(x,x)$ do not depend on the chart, so by Lemma~\ref{lem:sa:chart} we may use a
chart centered at $x_{0}$ with potential $\varphi=|z|^{2}+O(|z|^{3})$. By Lemma~\ref{lem:sa:real}, $\varphi$ is
realized near $0$ by a positive metric on $\mathcal O(1)\to\CP^{n}$, whose circle bundle $X'$ is a compact
regular Sasakian manifold with a BRT chart with the same potential (Remark~\ref{rem:sa:model}). By Proposition~\ref{tf:prop:H}(4) the
symbols of $X$ and $X'$ agree at $x_{0}$ modulo flat terms, and so does $|e(x_{0})|^{2}$
(Lemma~\ref{tf:lem:vol}). On $X'$, Proposition~\ref{tf:prop:diag} gives
$|e|^{2}\sigma_{j}(x,x)=a^{E}_{j}(\pi(x))/(2\pi^{n+1})$, with the coefficients of $\mathcal O(k)\otimes K_{\CP^{n}}$
for the metric realizing $\varphi$. These are the universal expressions evaluated on the K\"ahler metric
\cite[Thm.~4.1.2]{MM07}, which near $0$ is $g^{T}$.
\end{proof}

\begin{remark}[Consistency checks]\label{tf:rem:hanges}
(a) $a^{E}_{1}=\frac12r_{g}-r_{g}=-\frac12r_{g}$, since the curvature of $K_{M}$ is $-\mathrm{Ric}$. On $\CP^{n}$,
$E_{m}=\mathcal O(m-n-1)$ and $B^{E}_{m}=\pi^{-n}\prod_{i=1}^{n}(m-n-1+i)=(m/\pi)^{n}(1-\frac{n(n+1)}{2m}+O(m^{-2}))$,
so $a^{E}_{1}=-\frac{n(n+1)}{2}=-\frac12r_{g}$.

(b) Let $n=1$. For a real analytic rigid hypersurface $\operatorname{Im}w=f(z,\overline z)$ in $\mathbb C^{2}$,
Hanges \cite[Theorem]{Han90} computed the first three terms $h_{1},h_{0},h_{-1}$ of the symbol of the
Szeg\H o kernel, with kernels taken with respect to $dx\,dy\,dt$, where $z=x+\sqrt{-1}y$ and
$t=\operatorname{Re}w$ \cite[(2.2)]{Han90}. In BRT coordinates with $\vartheta=t$ and $\varphi=2f$, his phase is
the BRT phase, his volume form is $dV_{e}=d\lambda(z)\wedge d\vartheta$ (Lemma~\ref{tf:lem:vol}), and $u\mapsto u\,e$
identifies CR functions with CR $(1,0)$-forms. Since the symbol depends only on the local CR structure
and the volume form (\cite[Thm.~2.14, Prop.~3.19 and p.~162]{BouSj75}, as used in \cite[\S3]{Han90}), his
$h_{1-j}$ is, modulo flat terms, the term $\sigma_{j}$ of Proposition~\ref{tf:prop:H} for $n=1$. With
$g_{1\bar1}=\varphi_{z\bar z}=\det\varphi''$, and $r_{g}$ and $\Delta=g^{1\bar1}\partial_{z}\partial_{\bar z}$ as in
Lemma~\ref{lem:sa:n1}, his formulas \cite[(2.28), (2.32), (2.38)]{Han90} can be rewritten as
\[
 h_{1-j}(z,\overline z)=\frac{\det\varphi''}{2\pi^{2}}\,a^{E}_{j},\qquad j=0,1,2,\qquad
 a^{E}_{0}=1,\quad a^{E}_{1}=-\tfrac12 r_{g},\quad a^{E}_{2}=-\tfrac16\Delta r_{g}.
\]
These are the twisted coefficients for $n=1$: the formula of \cite[Rem.~1.8]{EXX} reduces to them,
because its terms quadratic in the curvature cancel when $n=1$. This agrees with Theorem~\ref{tf:thm:id}
for $j\le2$ and checks its normalizing constant independently.
\end{remark}

\begin{corollary}[The logarithmic term]\label{tf:cor:LT}
Let $W\subset X$ be open. The logarithmic term of $S$ vanishes to infinite order over $W$, in the sense of
Definition~\ref{cf:def:log} applied to the kernels $K$ in local frames, if and only if $a^{E}_{j}(g^{T})=0$ on
$W$ for every $j\ge n+1$.
\end{corollary}

\begin{proof}
By Proposition~\ref{tf:prop:H}, the proof of the equivalence of (i), (ii) and (iii) in
Theorem~\ref{thm:logS} applies to $K$, with Theorem~\ref{tf:thm:id} in place of Theorem~\ref{thm:introD}: it
uses the kernel only through \ref{K:rep}--\ref{K:sym}, and $|e|^{2}$ does not vanish. By
Remark~\ref{tf:rem:frame} the condition does not depend on the frame.
\end{proof}

\begin{proof}[Proof of Theorem~\ref{thm:introC}]
The self-adjointness of $A$ and the commutation with $S$ are Lemma~\ref{tf:lem:sa}. The representation of
$K$ is Proposition~\ref{tf:prop:H}, the statements on $\chi(k^{-1}A)S$ are Theorem~\ref{tf:thm:main}, the
identification of the $\sigma_{j}(x,x)$ is Theorem~\ref{tf:thm:id}, and the statement on the logarithmic
term is Corollary~\ref{tf:cor:LT}.
\end{proof}

\section{The regular case: Bergman kernels and ball quotients}\label{tf:sec:regular}
On a general complete Sasakian manifold there is no base and no line bundle, and $a^{E}_{j}(g^{T})$ is the
universal expression evaluated on the transverse metric, not the coefficient of a Bergman kernel. In this
section $X$ is the circle bundle of Section~\ref{tf:ssec:circle}. Then $\chi(k^{-1}A)S=\sum_{m}\chi(m/k)S_{m}$,
the $S_{m}$ are the Bergman projections of $E_{m}=L^{m}\otimes K_{M}$ (Lemma~\ref{tf:lem:fourier}), and the
$a^{E}_{j}$ are the coefficients of the Bergman functions $B^{E}_{m}$ (Proposition~\ref{tf:prop:diag}).

\begin{proposition}[Bergman and $(n,0)$ Szeg\H o]\label{tf:prop:BS}
Let $(L,h)\to M$ be as in Section~\ref{tf:ssec:circle}, $D\subset L^{*}$ the disc bundle, $X=\partial D$, and
$U\subset M$ open. Then:
\begin{enumerate}
\item the logarithmic term of the Szeg\H o projection on $(n,0)$-forms of $X$ vanishes to infinite order
      over $\pi^{-1}(U)$ if and only if $a^{E}_{j}=0$ on $U$ for all $j\ge n+1$;
\item the logarithmic coefficient of the Bergman kernel of $D$ on the diagonal vanishes to infinite
      order at $X$ over $U$ if and only if $a^{E}_{j}=0$ on $U$ for all $j\ge n+2$.
\end{enumerate}
So \textup{(1)} implies \textup{(2)}, and the converse holds exactly when in addition $a^{E}_{n+1}=0$ on $U$.
\end{proposition}

\begin{proof}
(1) is Corollary~\ref{tf:cor:LT} with Proposition~\ref{tf:prop:diag}. For (2), the Bergman kernel of $D$
on the diagonal is $e^{\varphi}\det g\,\frac1\pi\sum_{k\ge1}kB^{E}_{k}t^{k-1}$, $t=1-\rho$ \cite[Thm.~2.2]{EXX}.
This form requires neither compactness nor completeness: \cite[Thm.~2.2]{EXX} is stated for an arbitrary
polarized manifold, which by \cite[Def.~1.1]{EXX} need not be compact or complete; the decomposition of
$A^{2}_{(n+1,0)}(D)$ into fiberwise Fourier components is \cite[Prop.~2.3]{EXX}, the norm identity is
\cite[Lemma~2.4]{EXX}, and the convergence of the series on compact sets uses only the local bound
\cite[(2.7)]{EXX}. By Proposition~\ref{tf:prop:diag}, $kB^{E}_{k}\sim\pi^{-n}\sum_{j}a^{E}_{j}k^{n+1-j}$ locally
uniformly with derivatives, and the summation of Lemma~\ref{lem:sums} applies pointwise in $z$: the
indices $j\le n+1$ give poles in $\rho$, and $j=n+2+p$ gives $\rho^{p}\log\rho$ with a non-zero constant times
$a^{E}_{n+2+p}$. The difference between (1) and (2) is the term $j=n+1$: the Bergman series is the
derivative in $t$ of the Szeg\H o series $\sum_{k}B^{E}_{k}t^{k}$, and differentiation turns
$a^{E}_{n+1}\log\rho$ into a pole.
\end{proof}

Thus, in the regular case, the Bergman kernel of $D$ is governed by the same coefficients $a^{E}_{j}$ as
the Szeg\H o projection on $(n,0)$-forms of $X$, its logarithmic term starting at $j=n+2$ instead of
$j=n+1$, whereas the Szeg\H o projection on functions is governed by the $a^{L}_{j}$.

\begin{corollary}[Ball quotients]\label{tf:cor:ball}
Let $M=\Gamma\backslash\mathbb B^{n}$, with $\Gamma$ torsion free, let $(L,h)$ be a positive line bundle which is
locally the trivial bundle with weight $\varphi=-q\log(1-|z|^{2})$, $q>0$, and let
$\omega_{g}=\frac{\sqrt{-1}}{2}\partial\overline\partial\varphi$. Then the logarithmic terms of the Szeg\H o
projection on $(n,0)$-forms of the circle bundle $X$ and of the Bergman kernel of the disc bundle $D$
vanish to infinite order.
\end{corollary}

\begin{proof}
The Sasakian metric of $X$ is complete, since $M$ is complete and the fibers are circles. By
Corollary~\ref{tf:cor:LT} and Proposition~\ref{tf:prop:BS} it suffices to show $a^{E}_{j}=0$ for $j\ge n+1$.
These are universal local expressions in the metric and the curvatures \cite[Thm.~4.1.2]{MM07}, and $M$
is locally isometric to $\mathbb B^{n}$ with these data, so it suffices to compute on $\mathbb B^{n}$, where
Proposition~\ref{tf:prop:diag} applies. There $\det\varphi''=q^{n}(1-|z|^{2})^{-(n+1)}$, $|1|^{2}_{h}=(1-|z|^{2})^{q}$,
$|dz_{1}\wedge\dots\wedge dz_{n}|^{2}_{H}=(\det\varphi'')^{-1}$ and $dV_{g}=\det\varphi''\,d\lambda$, so the squared
$L^{2}$ norm of $f\,1^{m}\otimes dz_{1}\wedge\dots\wedge dz_{n}$ is $\int|f|^{2}(1-|z|^{2})^{qm}\,d\lambda$. Its
reproducing kernel on the diagonal is $\frac{\Gamma(n+qm+1)}{\pi^{n}\Gamma(qm+1)}(1-|z|^{2})^{-(n+1+qm)}$, and
multiplying by the pointwise squared norm $(1-|z|^{2})^{qm}q^{-n}(1-|z|^{2})^{n+1}$ of the frame gives
\[
 B^{E}_{m}=\frac{\Gamma(qm+n+1)}{(\pi q)^{n}\,\Gamma(qm+1)}=\pi^{-n}\prod_{i=1}^{n}\Bigl(m+\frac{i}{q}\Bigr),
\]
a polynomial of degree $n$ in $m$ with leading term $(m/\pi)^{n}$. Hence $a^{E}_{j}=0$ for $j\ge n+1$. As a
check, $a^{E}_{1}=n(n+1)/(2q)=-\frac12 r_{g}$, since $r_{g}=-n(n+1)/q$. The passage from the logarithmic
coefficient of the Bergman kernel of $D$ to Bergman logarithmic flatness of $X$, a local CR notion,
uses the localization \cite[Thm.~1.2, Cor.~1.4]{EXX}, which requires $M$ to admit a complete K\"ahler
metric, as it does here.
\end{proof}

\begin{remark}[K\"ahler--Einstein bases]\label{tf:lem:KE}
If the curvature of $K_{M}$ is $c$ times that of $L$ on an open set $U$, for a constant $c$, then on $U$,
$a^{E}_{j}=0$ for all $j\ge n+1$ if and only if $a^{L}_{j}=0$ for all $j\ge n+1$. Indeed, $L^{k}\otimes K_{M}$
then has curvature $(k+c)$ times that of $L$, and by the universality of the coefficients
\cite[Thm.~4.1.2]{MM07}, as formal series in $k^{-1}$,
$\sum_{j\ge0}a^{E}_{j}k^{n-j}=\sum_{j\ge0}a^{L}_{j}(k+c)^{n-j}$; the right-hand side is a polynomial in $k$ of
degree $n$ if and only if the left-hand side is. For the ball with $\varphi=-q\log(1-|z|^{2})$ this is the
relation between $B^{L}_{k}$ and $B^{E}_{k}$ in Remark~\ref{rem:models}, with $c=(n+1)/q$, and it gives a
second proof of Corollary~\ref{tf:cor:ball}.
\end{remark}

\subsection*{Acknowledgments}
This work was started at CY Cergy Paris Universit\'e, where the second author was a postdoctoral
researcher and the first author spent six months as a visitor. We thank Louis Ioos for his suggestions
and encouragement, which led to this work. The authors used
Claude (Anthropic) to assist in developing and checking
arguments and computations, in checking the literature, and in revising
the manuscript. The authors take full responsibility for the content of
this paper.

\appendix
\section{Proofs of the normal forms}\label{app:nf}
This appendix contains the proofs deferred from Section~\ref{sec:nf}, with two auxiliary lemmas and the
precise sources. We keep the notation, the inputs \ref{hyp:H1}--\ref{hyp:H2} and the hypotheses
\ref{K:rep}--\ref{K:sym} of that section.

\subsection{The BRT phase}\label{app:nf:phase}
\begin{remark}[The upper bound in (P3)]\label{cf:rem:upper}
Let $f$ be smooth and non-negative near a point of $\{z=w\}$, and let $f$ vanish on
$\{z=w\}$. Then $f\le C|z-w|^{2}$, locally uniformly in the remaining variables. Indeed every
point of $\{z=w\}$ is a minimum of $f$, so $\nabla f=0$ there, and Taylor's formula along the
segment from $(z,z)$ to $(z,w)$ gives
\[
 f(z,w)=\int_{0}^{1}(1-t)\,\partial_{t}^{2}f\bigl(z,z+t(w-z)\bigr)\,dt
 \le\tfrac12\sup|\nabla^{2}f|\;|z-w|^{2}.
\]
For $f=\delta$ this is the upper bound in (P3); the lower bound is strict pseudoconvexity.
The argument uses that $\delta$ vanishes on all of $\{z=w\}$, not only on the diagonal. For a
general phase $\psi$ this fails for $\operatorname{Im}\psi$, which need not vanish on $\{z=w\}$
where $\vartheta_{x}\ne\vartheta_{y}$.
\end{remark}

\begin{proof}[Proof of Proposition~\ref{nf:prop:phase}]
The BRT phase satisfies \textup{($\alpha$)--($\gamma$)} and the eikonal equation (P4), so by
\ref{hyp:H2} it is equivalent to $\psi$; by \cite[Def.~4.3, Thm.~4.4]{HM17JDG}, $\psi=g\phi+f$ with $g$
smooth and $f$ flat, and $g(x,x)=1$ since $d_{x}\psi(x,x)=d_{x}\phi(x,x)=\eta(x)$. The
substitution $t'=gt$ is justified by an almost analytic extension in $t$ and Stokes' formula, and the
flat term $f$ is removed after writing both phases in Malgrange form in the Reeb variable, as in
\cite[(3.1.50)--(3.1.58)]{CHH25}; this rewrites $I_{\psi}(\sigma)$ as $I(\sigma')$ with
$\sigma'_{j}=g^{-(n+1-j)}\sigma_{j}$.
 For circle bundles this is \cite[Thm.~3.8]{HM14}; for the Szeg\H o
projector of a compact Sasakian manifold see also \cite[(3.1.80)]{CHH25}.
\end{proof}

\subsection{Flat terms and uniqueness}\label{app:nf:normal}
\begin{lemma}[Taylor's formula at a complex point]\label{cf:lem:AA}
Let $f$ be smooth in a real variable $u$, with parameters, and let $\widetilde{f}$
be an almost analytic extension. For complex $c$ near the real axis,
$f(u)=\widetilde{f}(c)+(u-c)R(u)+E(u)$, where $R$ is smooth and $E$, with all
its derivatives, is $O(|\operatorname{Im}c|^{\infty})$.
\end{lemma}
\begin{proof}
Integrate $\frac{d}{ds}\widetilde{f}(c+s(u-c))$ over $[0,1]$. The
$\partial_{\zeta}$-part gives $(u-c)R$. The $\partial_{\bar\zeta}$-part is
$O(|\operatorname{Im}c|^{\infty})$, because
$\operatorname{Im}(c+s(u-c))=(1-s)\operatorname{Im}c$.
\end{proof}

\begin{proof}[Proof of Lemma~\ref{cf:lem:flat}]
For the phase of \cite[Prop.~1.1]{BouSj75}, (a) and (b) are contained in
\cite[(1.9), Rem.~1.10]{BouSj75}; the following direct argument applies to every Szeg\H o phase.
(a) If $a$ is flat at $z=w$,
which is the case in the applications below, a derivative of order $|\gamma|$ of
$e^{\sqrt{-1}t\psi}a$ is bounded by $C_{N}|z-w|^{N}t^{m+|\gamma|}e^{-ct|z-w|^{2}}$, whose
$t$-integral tends to $0$ as $z\to w$. In general, $|x-y|^{2}\le C|\psi|$ near the diagonal, by
Taylor's formula at $(y,y)$, since $d_{x}\psi(y,y)=\eta(y)$ and $\eta(T)=1$, together with
$|z-w|^{2}\le c^{-1}\operatorname{Im}\psi$; writing
$e^{\sqrt{-1}t\psi}=(\sqrt{-1}\psi)^{-M}\partial_{t}^{M}e^{\sqrt{-1}t\psi}$ and integrating by parts in
$t$, the flatness of $a$ absorbs $\psi^{-M}$ and the amplitude gains $t^{-M}$. (b) Put $D=\psi'-\psi$ and $\psi_{\tau}=\psi+\tau D$. Then
$I_{\psi'}(a)-I_{\psi}(a)=\int_{0}^{1}I_{\psi_{\tau}}(\sqrt{-1}tDa)\,d\tau$, and
$\operatorname{Im}\psi_{\tau}=(1-\tau)\operatorname{Im}\psi+\tau\operatorname{Im}\psi'\ge0$.
The amplitude $tDa$ is flat, so (a), applied uniformly in $\tau$, gives the
claim.
\end{proof}

\begin{lemma}[Flatness from rapid decay]\label{cf:lem:decay}
Let $B\in \mathscr{C}^{\infty}(U\times U)$ with $B(z,z)=0$ and $\operatorname{Im}B(z,w)\ge c|z-w|^{2}$,
and let $b\sim\sum_{j}b_{j}k^{m-j}$ be a classical symbol on $U\times U$,
depending smoothly on parameters. If $e^{\sqrt{-1}kB}b=O(k^{-\infty})$
uniformly, then every $b_{j}$ vanishes to infinite order at $z=w$.
\end{lemma}
\begin{proof}
Since $|e^{\sqrt{-1}kB}|=e^{-k\operatorname{Im}B}$ and, by Remark~\ref{cf:rem:upper},
$k\operatorname{Im}B(z,z+u/\sqrt k)$ is bounded for $u$ in compact sets, $b(z,z+u/\sqrt k,k)=O(k^{-\infty})$ there. Taylor
expansion of the $b_{j}$ at $w=z$ gives, for every $N$,
\[
 k^{-m}\,b(z,z+u/\sqrt k,k)=\sum_{s<N}k^{-s/2}
 \sum_{2j+|\alpha|=s}\frac{\partial^{\alpha}_{w}b_{j}(z,z)}{\alpha!}\,u^{\alpha}+O(k^{-N/2}),
\]
with $\alpha$ a multi-index in the real coordinates of $w$. Each polynomial in
$u$ on the right vanishes identically, and monomials of different degrees
are independent, so $\partial^{\alpha}_{w}b_{j}(z,z)=0$ for all $j$ and $\alpha$.
\end{proof}

\begin{proof}[Proof of Proposition~\ref{cf:prop:unique}]
Fix $(z,w,\vartheta_{y})$ near a diagonal point and put $\theta=\vartheta_{x}-\vartheta_{y}$.
Since $a$ does not depend on $\vartheta_{x}$, $I(a)$ is, as a tempered distribution in
$\theta\in\mathbb R$, the inverse Fourier transform of $t\mapsto2\pi\mathbf 1_{t>0}e^{\sqrt{-1}t\mathcal B}a(t)$:
\[
 I(a)=\int_{0}^{\infty}e^{\sqrt{-1}t(\theta+\mathcal B(z,w))}a(z,w,\vartheta_{y},t)\,dt,\qquad
 \int_{\mathbb R}e^{-\sqrt{-1}k\theta}I(a)\,d\theta=2\pi e^{\sqrt{-1}k\mathcal B(z,w)}a(z,w,\vartheta_{y},k)
\]
for $k>0$. Let
$\varrho\in \mathscr{C}^{\infty}_{c}(\mathbb R)$ equal $1$ near $0$, with support in the region where
$I(a)$ is smooth. On the support of $1-\varrho$ one has $\phi\ne0$, and integration
by parts in $t$ gives $\partial^{\ell}_{\theta}I(a)=O(|\theta|^{-2})$ for $|\theta|\ge1$, uniformly
near the diagonal; hence the Fourier transforms of $(1-\varrho)I(a)$ and of
$\varrho I(a)$ are $O(k^{-\infty})$. Therefore $e^{\sqrt{-1}k\mathcal B}a(k)=O(k^{-\infty})$ uniformly,
and Lemma~\ref{cf:lem:decay}, with $\operatorname{Im}\mathcal B\ge c|z-w|^{2}$ from (P3), gives the
claim. This is the semi-classical Fourier transform used in
\cite[\S8]{HM17JDG}; it relies on the phase being affine in $\vartheta_{x}$ and on $a$
not depending on $\vartheta_{x}$.
\end{proof}

\begin{proof}[Proof of Proposition~\ref{cf:prop:normal}]
\emph{Step 1: elimination of $\vartheta_{x}$.} Let $b\,t^{k}$ be a term of
$\sigma$. Lemma~\ref{cf:lem:AA} in $u=\vartheta_{x}$ at $c$ gives $b=\widehat{b}+\phi R+E$,
where $\widehat{b}(z,w,\vartheta_{y})=\widetilde{b}(z,c,w,\vartheta_{y})$ does not depend on
$\vartheta_{x}$, $R$ is smooth, and $E$ is flat by (P3). Since
$\phi e^{\sqrt{-1}t\phi}=-\sqrt{-1}\partial_{t}e^{\sqrt{-1}t\phi}$, one has
$I(\phi Rt^{k})\equiv I(\sqrt{-1}kRt^{k-1})$, the term with $\kappa'$ being smooth.
So a term of order $k$ becomes a term of $\mathcal C'$ of order $k$ plus a term
of order $k-1$, modulo a flat term, which is negligible by
Lemma~\ref{cf:lem:flat}; for $k=0$ the factor $k$ vanishes. Proceeding
through the orders and Borel-summing gives $\sigma'\in\mathcal C'$ with
$K\equiv I(\sigma')$, without terms of negative order if $\sigma$ has none.
This is the procedure of \cite[Lem.~3.3]{HS20}, in which the Malgrange
factorization is trivial.

\emph{Step 2: independence of $\vartheta_{y}$.} By \ref{K:inv},
$K(e^{uT}x,e^{uT}y)=K(x,y)$ for small $u$. By (P1), $e^{uT}$
shifts $\vartheta_{x}$ and $\vartheta_{y}$ by $u$ and preserves $\phi$, hence
$I(\sigma')(e^{uT}x,e^{uT}y)=I(\sigma'_{u})(x,y)$ with
$\sigma'_{u}(z,w,\vartheta_{y},t):=\sigma'(z,w,\vartheta_{y}+u,t)\in\mathcal C'$. Thus
$I(\sigma'_{u}-\sigma')\equiv0$ for $|u|$ small, and by
Proposition~\ref{cf:prop:unique} each
$\sigma'_{j}(z,w,\vartheta_{y}+u)-\sigma'_{j}(z,w,\vartheta_{y})$ is flat. So every Taylor
coefficient of $\sigma'_{j}$ at a diagonal point $(z_{0},z_{0})$ is locally
constant in $\vartheta_{y}$; hence $\partial_{\vartheta_{y}}\sigma'_{j}$ is flat, and so is
$\sigma'_{j}-\sigma'_{j}|_{\vartheta_{y}=0}$. The terms $\widehat{\sigma}_{j}:=\sigma'_{j}|_{\vartheta_{y}=0}$
define $\widehat{\sigma}\in\mathcal C$, and $K\equiv I(\widehat{\sigma})$ by
Lemma~\ref{cf:lem:flat}.
\end{proof}

\subsection{Almost sesquiholomorphy and polarization}\label{app:nf:CR}
As in Section~\ref{cf:sec:CR}, $K$ satisfies \ref{K:rep}--\ref{K:sym}, and $s\in\mathcal C$ is the symbol
given by Proposition~\ref{cf:prop:normal}.

\begin{proof}[Proof of Proposition~\ref{cf:prop:CR}]
Since $\overline{Z}_{l,x}K=0$ by \ref{K:CR}, and $\overline{Z}_{l}=\partial_{\bar z_{l}}$ on Reeb-independent
functions (\cite[(3.19)]{HS20}, exact here), (P4) gives
$0\equiv\overline{Z}_{l,x}I(s)=I(\partial_{\bar z_{l}}s)+I(t\varepsilon_{l}s)$. The last
term is smooth by Lemma~\ref{cf:lem:flat}, and
Proposition~\ref{cf:prop:unique} makes every $\partial_{\bar z_{l}}s_{j}$ flat.
Next, \ref{K:sym} and (P2) give $K\equiv I(s^{\dagger})$ with
$s^{\dagger}(z,w,t)=\overline{s(w,z,t)}\in\mathcal C$. Hence $s_{j}-s^{\dagger}_{j}$
is flat, and
$\partial_{w_{l}}s_{j}(z,w)=\overline{(\partial_{\bar z_{l}}s_{j})(w,z)}$
modulo flat functions.
\end{proof}

\begin{proof}[Proof of Corollary~\ref{cf:cor:polar}]
By Proposition~\ref{cf:prop:CR}, the Taylor series of $s_{j}$ at $(z_{0},z_{0})$
contains only monomials $(z-z_{0})^{\alpha}\overline{(w-z_{0})}^{\beta}$. These
monomials are independent on $w=z$, where the series becomes that of
$s_{j}(z,z)$. This is \cite[(3.21)--(3.29)]{HS20}, without the
$y_{2n+1}$-directions; for the amplitude of the semi-classical Toeplitz operator the same argument
is \cite[(3.1.85), (3.1.86)]{CHH25}.
\end{proof}

\begin{proof}[Proof of Proposition~\ref{nf:prop:symbol}]
Propositions~\ref{cf:prop:normal}, \ref{cf:prop:unique} and~\ref{cf:prop:CR} and
Corollary~\ref{cf:cor:polar}; for $\Pi$ the hypotheses hold, as explained after \ref{K:sym}. In two
charts the phases agree modulo flat terms (Lemma~\ref{lem:sa:chart}) and the symbols lie in the same
class, so by uniqueness they differ by flat terms; this gives the last assertion.
\end{proof}

\subsection{Sources}\label{app:nf:sources}
\begin{remark}[Sources]\label{nf:rem:sources}
The normal forms of Section~\ref{sec:nf} are known. For circle bundles, the phase is that of
\cite[Thm.~3.8]{HM14}, whose construction contains the Reeb-independent symbol \cite[\S3]{HM14}; the
uniqueness argument is the semi-classical Fourier transform along the Reeb flow of \cite[\S8]{HM17JDG};
for the semi-classical Toeplitz operator of a compact Sasakian manifold, the Reeb-independent
amplitude, its uniqueness and its determination by the diagonal values are
\cite[Thm.~1.3, Thm.~3.7, (3.1.85), (3.1.86)]{CHH25}; and the polarization is that of
\cite[(3.21)--(3.29)]{HS20}.
\end{remark}

\bibliographystyle{amsplain-abbrv}
\bibliography{LuTian_merged}

\providecommand{\bysame}{\leavevmode\hbox to3em{\hrulefill}\thinspace}
\providecommand{\MR}{\relax\ifhmode\unskip\space\fi MR }
\providecommand{\MRhref}[2]{%
  \href{http://www.ams.org/mathscinet-getitem?mr=#1}{#2}
}
\providecommand{\href}[2]{#2}
\begin{thebibliography}{10}

\bibitem{ALZ13}
C.~Arezzo, A.~Loi, and F.~Zuddas, \emph{Szeg{\"o} kernel, regular quantizations
  and spherical {CR}-structures}, Math. Z. \textbf{275} (2013), no.~3-4,
  1207--1216.

\bibitem{Bai57}
W.~L. Baily, \emph{On the imbedding of {$V$}-manifolds in projective space},
  Amer. J. Math. \textbf{79} (1957), 403--430.

\bibitem{BRT85}
M.~S. Baouendi, L.~P. Rothschild, and F.~Tr{\`e}ves, \emph{C{R} structures with
  group action and extendability of {CR} functions}, Invent. Math. \textbf{82}
  (1985), no.~2, 359--396. \MR{809720}

\bibitem{Bel01}
F.~A. Belgun, \emph{{Normal {CR} structures on compact 3-manifolds}}, Math. Z.
  \textbf{238} (2001), no.~3, 441--460. \MR{1869692}

\bibitem{BBS08}
R.~Berman, B.~Berndtsson, and J.~Sj\"{o}strand, \emph{A direct approach to
  {B}ergman kernel asymptotics for positive line bundles}, Ark. Mat.
  \textbf{46} (2008), no.~2, 197--217. \MR{2430724}

\bibitem{BouSj75}
L.~Boutet~de Monvel and J.~Sj\"{o}strand, \emph{Sur la singularit\'{e} des
  noyaux de {B}ergman et de {S}zeg{\H o}}, Journ\'{e}es \'{E}quations aux
  {D}\'{e}riv\'{e}es {P}artielles de {R}ennes (1975), Ast\'{e}risque, vol.
  34--35, Soc. Math. France, Paris, 1976, pp.~123--164. \MR{590106}

\bibitem{BoGa08}
C.~P. Boyer and K.~Galicki, \emph{Sasakian geometry}, Oxford Mathematical
  Monographs, Oxford University Press, Oxford, 2008. \MR{2382957}

\bibitem{Cat99}
D.~Catlin, \emph{The {B}ergman kernel and a theorem of {T}ian}, Analysis and
  geometry in several complex variables ({K}atata, 1997), Trends Math.,
  Birkh\"{a}user Boston, Boston, MA, 1999, pp.~1--23. \MR{1699887}

\bibitem{Cha26}
C.-C. Chang, \emph{{Quantization in Cauchy--Riemann geometry}}, PhD thesis,
  Universit\"at zu K\"oln, 2026.

\bibitem{CHH25}
C.-C. Chang, H.~Herrmann, and C.-Y. Hsiao, \emph{On the second coefficient in
  the semi-classical expansion of {T}oeplitz operators}, Anal. Math. Phys.
  \textbf{15} (2025), no.~4, Paper No. 112, 32. \MR{4938933}

\bibitem{CHT19}
J.-H. Cheng, C.-Y. Hsiao, and I.-H. Tsai, \emph{{Heat kernel asymptotics, local
  index theorem and trace integrals for {C}auchy-{R}iemann manifolds with
  {$S^1$} action}}, M\'em. Soc. Math. Fr. (N.S.) (2019), no.~162, vi+139.
  \MR{3975495}

\bibitem{CHT25}
\bysame, \emph{{Remarks on heat kernel asymptotics and local index theorem}},
  Bull. Inst. Math. Acad. Sin. (N.S.) \textbf{20} (2025), no.~2, 159--167.
  \MR{5005354}

\bibitem{CS18}
T.~C. Collins and G.~Sz\'{e}kelyhidi, \emph{K-semistability for irregular
  {S}asakian manifolds}, J. Differential Geom. \textbf{109} (2018), no.~1,
  81--109. \MR{3798716}

\bibitem{DLM06}
X.~Dai, K.~Liu, and X.~Ma, \emph{On the asymptotic expansion of {B}ergman
  kernel}, J. Differential Geom. \textbf{72} (2006), no.~1, 1--41. \MR{2215454}

\bibitem{DLM11}
\bysame, \emph{A remark on weighted {B}ergman kernels on orbifolds}, Math. Res.
  Lett. \textbf{19} (2012), no.~1, 143--148. \MR{2923181}

\bibitem{EXX}
P.~Ebenfelt, M.~Xiao, and H.~Xu, \emph{{Bergman kernels over polarized K\"ahler
  manifolds, Bergman logarithmic flatness, and a question of Lu--Tian}},
  arXiv:2510.22169 (2025).

\bibitem{Ebe18}
P.~Ebenfelt, \emph{{The log term in the {B}ergman and {S}zeg\H o{} kernels in
  strictly pseudoconvex domains in {$\mathbb{C}^2$}}}, Doc. Math. \textbf{23}
  (2018), 1659--1676. \MR{3890958}

\bibitem{EZ10}
M.~Engli{\v{s}} and G.~Zhang, \emph{{Ramadanov conjecture and line bundles over
  compact {H}ermitian symmetric spaces}}, Math. Z. \textbf{264} (2010), no.~4,
  901--912. \MR{2593299}

\bibitem{Fe74}
C.~Fefferman, \emph{The {B}ergman kernel and biholomorphic mappings of
  pseudoconvex domains}, Invent. Math. \textbf{26} (1974), 1--65. \MR{350069}

\bibitem{Gei97}
H.~Geiges, \emph{{Normal contact structures on {$3$}-manifolds}}, Tohoku Math.
  J. (2) \textbf{49} (1997), no.~3, 415--422. \MR{1464186}

\bibitem{Han90}
N.~Hanges, \emph{{Explicit formulas for the {S}zeg\H o{} kernel for some
  domains in {${\bf C}^2$}}}, J. Funct. Anal. \textbf{88} (1990), no.~1,
  153--165. \MR{1033918}

\bibitem{HeHsLi17}
H.~Herrmann, C.-Y. Hsiao, and X.~Li, \emph{{Szeg\H o{} kernel expansion and
  equivariant embedding of {CR} manifolds with circle action}}, Ann. Global
  Anal. Geom. \textbf{52} (2017), no.~3, 313--340. \MR{3711603}

\bibitem{HeHsLi18}
\bysame, \emph{Szeg\"{o} kernel asymptotic expansion on strongly pseudoconvex
  {CR} manifolds with {$S^1$} action}, Internat. J. Math. \textbf{29} (2018),
  no.~9, 1850061, 35. \MR{3845397}

\bibitem{HHL20}
\bysame, \emph{Torus equivariant {S}zeg{\H o} kernel asymptotics on strongly
  pseudoconvex {CR} manifolds}, Acta Math. Vietnam. \textbf{45} (2020), no.~1,
  113--135. \MR{4081369}

\bibitem{HeHsLiR}
\bysame, \emph{{Szeg\H o{} kernels and equivariant embedding theorems for {CR}
  manifolds}}, Math. Res. Lett. \textbf{29} (2022), no.~1, 193--246.
  \MR{4477684}

\bibitem{HHMS23}
H.~Herrmann, C.-Y. Hsiao, G.~Marinescu, and W.-C. Shen, \emph{Semi-classical
  spectral asymptotics of {T}oeplitz operators on {CR} manifolds}, Mem. Eur.
  Math. Soc., to appear; arXiv:2303.17319. Theorem numbers refer to the
  accepted version: its Theorem~1.1.1, Corollary~1.1.2, (1.1.4) and
  Theorem~4.2.3 are Theorem~1.1, Corollary~1.2, (1.2) and Theorem~5.5 of
  arXiv:2303.17319v2, and Theorem~4.3.1 is not in v2.

\bibitem{Hir93}
K.~Hirachi, \emph{{Scalar pseudo-{H}ermitian invariants and the {S}zeg\H o{}
  kernel on three-dimensional {CR} manifolds}}, Complex geometry ({O}saka,
  1990), Lecture Notes in Pure and Appl. Math., vol. 143, Dekker, New York,
  1993, pp.~67--76. \MR{1201602}

\bibitem{Hoe03}
L.~H\"{o}rmander, \emph{The analysis of linear partial differential operators.
  {I}}, Classics in Mathematics, Springer-Verlag, Berlin, 2003, Distribution
  theory and Fourier analysis, Reprint of the second (1990) edition.
  \MR{1996773}

\bibitem{HLM23}
C.-Y. Hsiao, X.~Li, and G.~Marinescu, \emph{{Semi-classical asymptotics of
  partial Bergman kernels on $\mathbb R$-symmetric complex manifolds with
  boundary}}, arXiv:2208.12412.

\bibitem{Hs10}
C.-Y. Hsiao, \emph{Projections in several complex variables}, M\'{e}m. Soc.
  Math. Fr. (N.S.) (2010), no.~123, 131. \MR{2780123}

\bibitem{HLM21}
C.-Y. Hsiao, X.~Li, and G.~Marinescu, \emph{Equivariant {K}odaira embedding for
  {CR} manifolds with circle action}, Michigan Math. J. \textbf{70} (2021),
  no.~1, 55--113. \MR{4255089}

\bibitem{HM14}
C.-Y. Hsiao and G.~Marinescu, \emph{Asymptotics of spectral function of lower
  energy forms and {B}ergman kernel of semi-positive and big line bundles},
  Comm. Anal. Geom. \textbf{22} (2014), no.~1, 1--108. \MR{3194375}

\bibitem{HM17JDG}
\bysame, \emph{On the singularities of the {S}zeg{\H o} projections on lower
  energy forms}, J. Differential Geom. \textbf{107} (2017), no.~1, 83--155.
  \MR{3698235}

\bibitem{HMW22}
C.-Y. Hsiao, G.~Marinescu, and H.~Wang, \emph{{Szeg\H o kernel asymptotics on
  complete strictly pseudoconvex {CR} manifolds}}, J. Geom. Anal. \textbf{32}
  (2022), no.~11, Paper No. 266, 53 pp.

\bibitem{HS20}
C.-Y. Hsiao and W.-C. Shen, \emph{On the second coefficient of the asymptotic
  expansion of {B}outet de {M}onvel--{S}j\"{o}strand}, Bull. Inst. Math. Acad.
  Sin. (N.S.) \textbf{15} (2020), no.~4, 329--365. \MR{4205414}

\bibitem{HL}
X.~Huang and X.~Li, \emph{{Bergman-{E}instein metric on a {S}tein space with a
  strongly pseudoconvex boundary}}, Comm. Anal. Geom. \textbf{31} (2023),
  no.~7, 1669--1692. \MR{4785570}

\bibitem{Kaw79}
T.~Kawasaki, \emph{The {R}iemann--{R}och theorem for complex {$V$}-manifolds},
  Osaka J. Math. \textbf{16} (1979), no.~1, 151--159.

\bibitem{Koh86}
J.~J. Kohn, \emph{The range of the tangential {C}auchy-{R}iemann operator},
  Duke Math. J. \textbf{53} (1986), no.~2, 525--545. \MR{850548}

\bibitem{Lem92}
L.~Lempert, \emph{On three-dimensional {C}auchy-{R}iemann manifolds}, J. Amer.
  Math. Soc. \textbf{5} (1992), no.~4, 923--969. \MR{1157290}

\bibitem{LUZ16}
A.~Loi, D.~Uccheddu, and M.~Zedda, \emph{On the {S}zeg{\"o} kernel of
  {C}artan--{H}artogs domains}, Ark. Mat. \textbf{54} (2016), no.~2, 473--484.

\bibitem{LMZ19}
A.~Loi, R.~Mossa, and F.~Zuddas, \emph{Finite {TYCZ} expansions and c{s}c{K}
  metrics}, J. Math. Anal. Appl. \textbf{484} (2020), no.~1, Paper No. 123715,
  20, arXiv:1903.07679. \MR{4037593}

\bibitem{LPZ26}
A.~Loi, M.~Palmieri, and F.~Zuddas, \emph{{The Szeg\H o log-term and local
  rigidity of K\"ahler--Einstein polarizations of compact Hermitian symmetric
  spaces}}, 2026, arXiv:2611.12708.

\bibitem{Lu00}
Z.~Lu, \emph{On the lower order terms of the asymptotic expansion of
  {T}ian-{Y}au-{Z}elditch}, Amer. J. Math. \textbf{122} (2000), no.~2,
  235--273. \MR{1749048}

\bibitem{LS26}
Z.~Lu and Y.~Shi, \emph{{Local rigidity of Bergman logarithmic flatness on
  compact Hermitian symmetric spaces}}, 2026, arXiv:2611.14507.

\bibitem{LT}
Z.~Lu and G.~Tian, \emph{{The log term of the {S}zeg\H o{} kernel}}, Duke Math.
  J. \textbf{125} (2004), no.~2, 351--387. \MR{2096677}

\bibitem{MM07}
X.~Ma and G.~Marinescu, \emph{Holomorphic {M}orse inequalities and {B}ergman
  kernels}, Progress in Mathematics, vol. 254, Birkh\"{a}user Verlag, Basel,
  2007. \MR{2339952}

\bibitem{MSY08}
D.~Martelli, J.~Sparks, and S.-T. Yau, \emph{{S}asaki--{E}instein manifolds and
  volume minimisation}, Comm. Math. Phys. \textbf{280} (2008), no.~3, 611--673.
  \MR{2399609}

\bibitem{MeSj74}
A.~Melin and J.~Sj\"{o}strand, \emph{Fourier integral operators with
  complex-valued phase functions}, Fourier integral operators and partial
  differential equations ({C}olloq. {I}nternat., {U}niv. {N}ice, {N}ice, 1974),
  Lecture Notes in Math., Vol. 459, Springer, Berlin, 1975, pp.~120--223.
  \MR{0431289}

\bibitem{MenSj78}
A.~Menikoff and J.~Sj\"{o}strand, \emph{On the eigenvalues of a class of
  hypoelliptic operators}, Math. Ann. \textbf{235} (1978), no.~1, 55--85.
  \MR{481627}

\bibitem{NN57}
A.~Newlander and L.~Nirenberg, \emph{Complex analytic coordinates in almost
  complex manifolds}, Ann. of Math. (2) \textbf{65} (1957), 391--404.

\bibitem{RT11_1}
J.~Ross and R.~Thomas, \emph{Weighted {B}ergman kernels on orbifolds}, J.
  Differential Geom. \textbf{88} (2011), no.~1, 87--107. \MR{3019812}

\bibitem{Ruk95}
P.~Rukimbira, \emph{Chern--{H}amilton's conjecture and {$K$}-contactness},
  Houston J. Math. \textbf{21} (1995), no.~4, 709--718.

\bibitem{Sco83}
P.~Scott, \emph{{The geometries of {$3$}-manifolds}}, Bull. London Math. Soc.
  \textbf{15} (1983), no.~5, 401--487. \MR{705527}

\bibitem{Tan69}
S.~Tanno, \emph{{Sasakian manifolds with constant {$\phi$}-holomorphic
  sectional curvature}}, Tohoku Math. J. (2) \textbf{21} (1969), 501--507.
  \MR{251667}

\bibitem{Ti90}
G.~Tian, \emph{On a set of polarized {K}\"{a}hler metrics on algebraic
  manifolds}, J. Differential Geom. \textbf{32} (1990), no.~1, 99--130.
  \MR{1064867}

\bibitem{Zel98}
S.~Zelditch, \emph{Szeg{\H o} kernels and a theorem of {T}ian}, Internat. Math.
  Res. Notices (1998), no.~6, 317--331. \MR{1616718}

\end{thebibliography}

\end{document}